\documentclass[11pt]{amsart}
\usepackage{amsmath,amsfonts,amssymb,amsthm,amsbsy,bbm}
\usepackage[curve]{xypic}
\usepackage{longtable}
\usepackage{array}
\newtheorem{theorem}[equation]{Theorem}
\newtheorem{proposition}[equation]{Proposition}

\newtheorem{corollary}[equation]{Corollary}
\newtheorem{lemma}[equation]{Lemma}
\newtheorem{remark}[equation]{Remark}
\numberwithin{equation}{section}

\newcommand\ii{\mathrm{i}}
\newcommand\into{\hookrightarrow}

\newcommand\bla{{\boldsymbol{\lambda}}}
\newcommand\centl{{C_\bla}}
\newcommand\om{^{(\omega)}}
\newcommand\pom{^{\prime(\omega)}}
\newcommand\trivialrep{{1_G}}
\newcommand\Id{{\mathrm{Id}}}
\newcommand\reg{{\mathrm{reg}}}

\DeclareMathOperator{\Ad}{\mathrm{Ad}}
\DeclareMathOperator{\Aut}{\mathrm{Aut}}
\DeclareMathOperator{\End}{\mathrm{End}}
\DeclareMathOperator{\Br}{\mathrm{Br}}
\DeclareMathOperator{\Diag}{\mathrm{Diag}}
\DeclareMathOperator{\Gal}{\mathrm{Gal}}
\DeclareMathOperator{\Hom}{\mathrm{Hom}}
\DeclareMathOperator{\Image}{\mathrm{Im}}
\DeclareMathOperator{\Irr}{\mathrm{Irr}}
\DeclareMathOperator{\Ker}{\mathrm{Ker}}
\DeclareMathOperator{\Out}{\mathrm{Out}}
\DeclareMathOperator{\Spec}{\mathrm{Spec}}
\DeclareMathOperator{\Trace}{\mathrm{Trace}}
\DeclareMathOperator{\Ind}{\mathrm{Ind}}
\DeclareMathOperator{\Res}{\mathrm{Res}}
\DeclareMathOperator{\GL}{\mathrm{GL}}
\DeclareMathOperator{\PGL}{\mathrm{PGL}}
\DeclareMathOperator{\Sp}{\mathrm{Sp}}
\DeclareMathOperator{\SO}{\mathrm{SO}}
\DeclareMathOperator{\sgn}{\mathrm{sgn}}
\DeclareMathOperator{\can}{\mathrm{Nab}}
\newcommand{\odd}{{2'}}
\newcommand{\BC}{{\mathbbm C}}
\newcommand{\BF}{{\mathbbm F}}
\newcommand{\BQ}{{\mathbbm Q}}
\newcommand{\BZ}{{\mathbbm Z}}
\newcommand{\BR}{{\mathbbm R}}
\newcommand{\BN}{{\mathbbm N}}
\newcommand{\CE}{{\mathcal E}}
\newcommand{\CH}{{\mathcal H}}
\newcommand{\CL}{{\mathcal L}}
\newcommand{\CS}{{\mathcal S}}
\newcommand{\CT}{{\mathcal T}}
\newcommand{\bT}{{\mathbf T}}
\newcommand{\bs}{{\mathbf s}}
\newcommand{\bt}{{\mathbf t}}
\newcommand{\bu}{{\mathbf u}}
\newcommand{\Abar}{{\overline A}}
\newcommand{\Cbar}{{\overline C}}
\newcommand{\Nbar}{{\overline N}}
\newcommand{\inv}{^{-1}}
\def\lexp#1#2{\kern\scriptspace\vphantom{#2}^{#1}\kern-\scriptspace#2}
\newcommand{\eg}{{\it e.g}}
\newcommand{\cf}{{\it cf}}
\newcommand{\ie}{{\it i.e}}
\newcommand{\cc}{{\mathfrak{c}}}
\newcommand{\Sgot}{{\mathfrak S}}
\newcommand{\Agot}{{\mathfrak A}}
\def\gal#1{{\mathfrak s_{#1}}}
\newcommand{\GAP}{{\tt GAP}}
\newcommand{\CHEVIE}{{\tt CHEVIE}}
\newcommand{\MAGMA}{{\tt MAGMA}}

\newcommand{\ibar}{{\overline\iota}}
\newcommand{\itilde}{{\tilde\iota}}

\begin{document}
\author{I.~Marin \and J.~Michel}
\address{Institut  de math{\'e}matiques de Jussieu, universit{\'e} Paris VII, 175,
rue du Chevaleret, 75013 Paris}
\email{marin@math.jussieu.fr}
\email{jmichel@math.jussieu.fr}
\title[Galois  automorphisms]{Automorphisms  of complex reflection groups}
\date{March 7, 2009}
\subjclass[2000]{20F55, 20F28, 20C99}
\begin{abstract}
Let  $G\subset\GL(\BC^r)$ be a  finite complex reflection  group. We show that
when  $G$ is irreducible, apart from the exception $G=\Sgot_6$, as well as for
a  large  class  of  non-irreducible  groups,  any  automorphism of $G$ is the
product  of a central automorphism and of an automorphism which preserves the
reflections.  We  show  further  that  an  automorphism  which  preserves  the
reflections  is  the  product  of  an  element of $N_{\GL(\BC^r)}(G)$ and of a
``Galois'' automorphism: we show that $\Gal(K/\BQ)$, where $K$ is the field of
definition  of  $G$, injects into  the group  of outer  automorphisms of $G$,
and that this injection can be chosen such that it induces
the usual Galois action on characters of $G$, apart from a few exceptional
characters; further, replacing if needed $K$ by an
extension  of degree $2$, the  injection can be lifted  to $\Aut(G)$, and every
irreducible representation admits a model which is equivariant with respect to
this lifting. Along the way we show that the fundamental invariants of $G$ can
be chosen rational.
\end{abstract}
\maketitle
\section{Introduction}

Let  $G$ be a  finite complex reflection  group of rank  $r$, that is a finite
group  generated by  (pseudo)reflections in  a vector  space $V$  of dimension
$r<\infty$  over  $\BC$.  We  say  that  $G$  is  {\em  irreducible}  when $G$
acts irreducibly on $V$.  In this  paper, we  determine the  group $\Aut(G)$ of
automorphisms  of $G$ when $G$ is irreducible, as well as for a large class of
non-irreducible groups.

An automorphism $\alpha\in\Aut(G)$ such that $\alpha(g)\in g ZG$ for all $g\in
G$, where $ZG$ denotes the center of $G$,  is called  {\em central}.  
These automorphisms  form a  normal subgroup of $\Aut(G)$, denoted $C$.

We call {\em reflections} of $G$ all the elements which act as a
(pseudo)-reflection on $V$ and denote by $A$ the subgroup of automorphisms of
$G$ which preserve the reflections. One of our main results is that
\begin{theorem}\label{auts}
Let $G$ be an irreducible complex reflection group different from
$G(1,1,6)\simeq \Sgot_6$; then $\Aut(G)=C\cdot A$.
If $\dim V>2$, then $\Aut(G)=C\rtimes A$.
\end{theorem}
Further, we extend this result to a large class of non-irreducible groups (see
\ref{KRS  for quasi-indec}).

For  $G$ of rank 2, it may happen that $C$ intersects non-trivially
$A$, and for $G$ of rank 1 we have $C=A$.

The  main observation of this paper is the interpretation of $A$
as ``Galois'' automorphisms.

Let $\chi_V$ be the character of $G$ afforded by $V$.
It  is well known (see \eg.  \cite[7.1.1]{benson}) that the representation $V$
of  $G$ can be realized over the subfield $K\subset \BC$ generated by
the  values of $\chi_V$ ; the  field $K$ is called the {\it field of
definition}  of $G$. Further, by  a well known  result of Benard and Bessis
(\cite{benard},  \cite{bessis}) {\it any}  irreducible representation of
$G$ can be realized over $K$. 

Let  $\Gamma=\Gal(K/\BQ)$,  and  let  $\gamma\in\Gamma$.  If  $a_\gamma$ is an
automorphism  of $G$ such that for any  $g\in G$ we have $\chi_V\circ a_\gamma
=\gamma(\chi_V)$  we  say  that  $a_\gamma$  is  a  {\em  Galois automorphism}
corresponding  to $\gamma$. 

\begin{theorem} \label{thref}
Assume $G$ irreducible. Then any faithful irreducible representation where
the reflections of $G$ still act by reflections is a Galois conjugate of $V$.
\end{theorem}

As a corollary, for an irreducible group, $A$ identifies to the group of
Galois automorphisms.

Note that $a_\gamma$  is determined by $\gamma$ up
to  an  element  of  the  normal  subgroup  $N$  formed of the Galois
automorphisms corresponding to the identity element of $\Gamma$. Equivalently,
$N$ is the subgroup of automorphisms induced by an element of $N_{\GL(V)}(G)$;
they  contain  the  inner  automorphisms.  We  call  the  induced outer
automorphisms {\em
diagram}  automorphisms (since they extend the corresponding notion for finite
Coxeter groups).

What is remarkable is that any element of the Galois group can be represented
by an automorphism. 
We  denote  by  $\Out(G)$  the  group  of outer automorphisms, identified with
$\Aut(G)/(G/ZG)$.  
Let $\Abar$ (resp. $\Nbar$) be the image of 
$A$ (resp. $N$) in $\Out(G)$.

\begin{theorem} \label{A->Gal surjective}
Assume $G$ irreducible. The natural map $\Abar\to\Gamma$ is surjective
and admits a section, thus $\Abar\simeq \Nbar \rtimes \Gamma$.
\end{theorem}

Theorem \ref{A->Gal surjective} can be reformulated as an equivariance
property of $\chi_V$:
there exists an injection $\ibar_V:\Gamma\to\Out(G)$ such that 
$\chi_V\circ\ibar_V(\gamma)=\gamma(\chi_V)$.

It is actually a special case of an equivariance property that
can be proved for almost any character. We call a representation
\emph{rational} if it admits a (matrix) model over $\BQ$.

\begin{theorem}\label{A}   Assume   $G$   irreducible.   For  any  irreducible
character  $\chi$ of $G$  but eight exceptions  there exists an injective
homomorphism:  $\Gamma  \xrightarrow{\ibar_\chi}\Out(G)$  such  that  for  any
$\gamma\in\Gamma$,  we have $\chi\circ\ibar_\chi(\gamma)=\gamma(\chi)$.

The  exceptions  are  the  characters  of  the  four  (rational) 5-dimensional
representations  of  $G_{27}$,  two  of  the  3  rational  representations  of
dimension 6 of $G_{29}$, and the two rational representations of dimension 120
of  $G_{34}$. In these cases,  the complex conjugation $\cc$  is an element of
$\Gamma$,  and  its  image  in  any injective morphism $\Gamma\to\Out(G)$ will
exchange in pairs the mentioned characters.

For an exceptional group (and a non-exceptional representation), or an
imprimitive group $G(de,e,r)$ where $\gcd(e,r)\le 2$, there is a unique way
to take the injection
$\ibar_\chi$ independent of $\chi$; we will call $\ibar$ the common value. 
In the remaining cases, we can make choices of $\ibar_\chi$ which are in a 
single orbit for $\Nbar$-conjugacy.
\end{theorem}

The  numbers  above  refer  to  the  Shephard-Todd classification of irreducible
complex  reflection groups. 

Theorem \ref{A} has the following consequence, which can be formulated without
exception and without assuming $G$ irreducible.

\begin{corollary}\label{critere  rationalite}  Let  $G$  be any finite complex
reflection  group, and let $\chi$ be an irreducible character of $G$. If,
for all $a \in \Out(G)$ we have $\chi \circ a =\chi$, then $\chi$ takes its
values in $\BQ$.
\end{corollary}

The  map $\ibar$ can  ``almost'' be lifted to
$\Aut(G)$. More precisely we have

\begin{theorem}\label{xx}  For $G$ irreducible and $\chi$ outside the
exceptions of theorem \ref{A}, there
exists  an extension  $K'$ of  $K$, which  is abelian  over $\BQ$  and at most
quadratic over $K$, and an injective homomorphism
$\itilde_\chi:\Gal(K'/\BQ)\to\Aut(G)$    such   that    the   composition   of
$\itilde_\chi$   with  the  natural  epimorphism  $\Aut(G)\to\Out(G)$  factors
through $\Gamma$ and induces $\ibar_\chi$.

One  may take $K'=K$  except when $G$  is $G_{27}$ or  $G$ is a dihedral group
$G(e,e,2)$  such that $\BQ(\zeta_e)$ contains no quadratic imaginary extension
of $\BQ$.

The  map  $\itilde_\chi$  can  be  taken  independent  of  $\chi$  in the same
circumstances  as $\ibar_\chi$; in such a case we will write $\itilde$ for the
common value.
\end{theorem}

Suppose  we can find a model over  some field $L$ of a faithful representation
$\rho$  of  $G$  such  that  $\rho(G)$  as  a  set  is  globally  invariant by
$\Gal(L/\BQ)$.  Then  any  $\gamma\in\Gal(L/\BQ)$  induces  a  permutation  of
$\rho(G)$ and thus a permutation on $G$, which is an automorphism. We get thus
an  homomorphism $\itilde:\Gal(L/\BQ)\to\Aut(G)$. If  $L=K'$ this $\itilde$
is equivariant with respect to the character of $\rho$ thus is suitable for
the $\itilde$ of theorem \ref{xx}.

Given  an homomorphism $\Gal(L/\BQ)\xrightarrow\itilde\Aut(G)$, and a model of
an  arbitrary  representation  $\rho$  over  $L$,  we  say  that  the model is
$\itilde$-equivariant  if for any $\gamma\in\Gal(L/\BQ)$  and any $g\in G$, we
have  $\rho(\itilde(\gamma)(g))=\gamma(\rho(g))$. Note that
this implies the global invariance of $\rho(G)$ under $\Gal(L/\BQ)$.
Considerations  as above led us to check the following theorem :

\begin{theorem}\label{C}  Let $G$ be an  irreducible finite complex reflection
group  which is not $G_{22}$ and let  $K'$ and $\itilde_{\chi}$ be as in theorem
\ref{xx}.  Let $\rho$ be any irreducible representation of $G$, with character
$\chi$, which is not an
exception  in \ref{A}.  Then there  is a  model of  $\rho$ over  $K'$ which is
$\itilde_{\chi}$-equivariant. \end{theorem}

In the case of $G_{22}$, we have $K=K'=\BQ(\ii,\sqrt 5)$. There is no globally
invariant model over $K'$ of the reflection representation. We need to replace
$K'$ by the extension $K''=\BQ(e^{\frac{2i\pi}{20}})$ to get an
invariant   model  (then  $K''$  also  works  for  any  other
representation).   This  example  illustrates  the  fact  that  asking  for  a
$\itilde_{\chi}$-equivariant  model is stronger than merely asking as in theorem
\ref{xx}  for a  model such  that $\chi\circ\itilde_\chi(\gamma)=\gamma(\chi)$
where $\chi$ is the character of $\rho$.

Theorem \ref{C} has the following consequence (which could have been
observed from the explicit values of the invariants discovered by various
authors, mostly in the 19th century).
In this proposition we take $K''$ as above for $G_{22}$, and in the other
cases we let $K'' = K'$.

\begin{corollary}\label{invariants}  Let $G\subset\GL(V)$ be any irreducible
complex  reflection group where $V$ is a $K''$-vector space, where $K''$ is as
above.  There is $\BQ$-form  $V=V_0\otimes_\BQ K''$ such  that the fundamental
invariants  of $G$  can be  taken rational,  \ie. in  the symmetric algebra of
the dual of $V_0$.
\end{corollary}

Let $V^\reg$ be the complement in $V$ of the reflecting hyperplanes for $G$.
The fundamental group of the variety $V^\reg/G$ is the {\em braid group} of
$G$. Since the variety $V/G$ is an affine
space it is defined over $\BQ$. A geometric reformulation
of \ref{invariants} is that the morphism $V \to V/G$ is also defined over
$\BQ$. We deduce from \ref{invariants} the following
\begin{corollary}\label{quotient} The varieties $V^\reg$, $V^\reg/G$
as well as the quotient morphism $V^\reg\to V^\reg/G$ are defined
over $\BQ$.
\end{corollary}

The  techniques we use for  checking theorem \ref{C} enable  us also to give a
probably   shorter   proof   of   the   Benard-Bessis   theorem   (see  remark
\ref{benardbessis}).

Finally  we  extend  our  results  to  a large class of non-irreducible finite
complex  reflection groups. We say that a group $G$ is {\em indecomposable} if
it   admits  no  non-trivial  decomposition  as  a  direct  product.  We  show
\begin{proposition}\label{quasi-indec} An irreducible complex reflection group
has  a decomposition $G=Z\times \hat{G}$ where $Z\subseteq ZG$, and where $\hat{G}$
is indecomposable non-abelian or trivial. In this decomposition $Z$ is unique.
\end{proposition}

We call $\hat{G}$ the non-abelian factor of $G$, denoted $\can(G)$; it is unique up
to isomorphism, since $Z$ is unique.
We then show:
\begin{theorem}\label{KRS  for quasi-indec} Let $G=G_1\times\ldots G_n$ be the
decomposition   in  irreducible   factors  of   a  complex   reflection  group
$G\subset\GL(V)$.  Then the following are equivalent:
\begin{itemize}
\item Any  automorphism of  $G$ is  the product of a central
automorphism  and of an  automorphism which preserves  the reflections.
\item No $G_i$  is isomorphic to  $\Sgot_6$, and the  decomposition has the
property  that for  any $i,j$,  we have  $ \can(G_i)\simeq\can(G_j)\Rightarrow
G_i\simeq G_j$.
\end{itemize}
\end{theorem}

We give the list of central and non-abelian factors of irreducible groups in
table \ref{bessisaut}.
\medskip

\subsection*{Acknowledgements}

We  thank David Harari for teaching us about the Brauer group,
Gunter Malle for a thorough reading of a first version of the paper, and
Luis Paris for suggesting the use of the Krull-Remak-Schmidt theorem.

\section{\label{section notation}Background on the infinite series}
We  will use  the following  notation throughout  this paper.  $V$ is a vector
space  over the subfield  $K$ of $\BC$.  $G$ is a  finite subgroup of $\GL(V)$
generated  by (pseudo)-reflections.  We assume  that $K$  is the {\sl field of
definition}  of $G$,  that is,  it is  the subfield  of $\BC$ generated by the
traces  of the  elements of  $G$. The  group $G$  is {\sl  irreducible} if the
representation  $V$  of  $G$  is  irreducible.  

We  denote  by  $g\mapsto\Ad  g$  the  map  which  maps  $g\in G$ to the inner
automorphism $x\mapsto \lexp gx=gxg\inv$.

The  Shephard-Todd classification of irreducible finite reflection groups (see
\eg.  \cite{cohen}) shows  that there  is one  infinite series,  identified by
three  positive integer parameters $d,e$ and  $r$, and denoted by $G(de,e,r)$.
In  addition  there  are  34  exceptional  
ones,  denoted  $G_4$ to $G_{37}$.

For the convenience of the reader, and for future reference,
we review the definition and some properties of the groups $G(de,e,r)$, 

\subsection{\label{Gde,e,r} The groups $G(de,e,r)$}

We  will  denote  by  $\mu_n$  the  subgroup  of  $n$-th  roots  of  unity  in
$\BC^\times$,and by $\zeta_n$ the primitive root of unity $e^{2i\pi/n}$.

The group $G(de,e,r)$ is defined as the subgroup of $\GL_r(\BC)$ consisting of
$r$  by $r$ monomial  matrices with entries  in $\mu_{de}$, and  such that the
product  of  non-zero  entries  is  in  $\mu_d$;  this  defines an irreducible
reflection  group except for the case  of $\Sgot_r=G(1,1,r)$, where one has to
quotient  by  the  one-dimensional  fixed  points,  and  for  $G(2,2,2) \simeq
\Sgot_2^2$.  


The  field of  definition $K$  of $G(de,e,r)$  is the  field generated  by the
traces  of elements of $G(de,e,r)$. We get $K=\BQ(\zeta_{de})$, except when $d
= 1$ and $r = 2$. In this last case, $G(e,e,2)$ is the dihedral group of order
$2e$,  whose  field  of  definition  is  $K=\BQ(\cos(\frac{2  \pi}{e})) = \BQ(
\zeta_e + \zeta_e\inv)$.

The  group $\Sgot_r=G(1,1,r)$  is a  subgroup of  $G(de,e,r)$; if we denote by
$D(de,e,r)$  the  subgroup  of  diagonal  matrices  in $G(de,e,r)$, we have a
semi-direct product decomposition $G(de,e,r)=D(de,e,r)\rtimes\Sgot_r$.

The   group  $G(d,1,r)$  is   generated  by  the   set  of  reflections  $\{
t,s_1,\ldots,s_{r-1}\}$ where $t$ is the matrix $\Diag(\zeta_d,1,\ldots,1) \in
D(d,1,r)$  and $s_k\in\Sgot_r$ is the  permutation matrix corresponding to the
transposition $(k,k+1)$.

The  group $G(de,e,r)$ is a normal subgroup of index $e$ in $G(de,1,r)$; it is
generated  by  the  set  of  reflections $\{t',s'_1,s_1,\ldots,s_{r-1}\}$ where
$t'=t^e$    and   $s'_1=s_1^t$;    here   $t$    refers   to   the   generator
$t=\Diag(\zeta_{de},1,\ldots,1)$  of $G(de,1,r)$. The generator $t'$ should be
omitted  (being  trivial)  if  $d=1$  and  $s'_1$  should  be  omitted  (being
unnecessary) if $e=1$. The quotient $G(de,1,r)/G(de,e,r)$ is cyclic, generated
by the image of $t$.

\subsection{Presentations}

The following diagram gives a presentation of $G(d,1,r)$:
\def\ncnode#1#2{{\kern -0.4pt\mathop\bigcirc\limits_{#2}\kern-8.6pt{\scriptstyle#1}\kern 2.3pt}}
\def\dbar#1pt{{\rlap{\vrule width#1pt height2pt depth-1pt} 
                 \vrule width#1pt height4pt depth-3pt}}
\def\nnode#1{{\kern -0.6pt\mathop\bigcirc\limits_{#1}\kern -1pt}}
\def\bar#1pt{{\vrule width#1pt height3pt depth-2pt}}
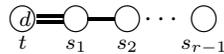
$$\ncnode dt\dbar10pt\nnode{s_1}\bar10pt\nnode{s_2}\cdots\nnode{s_{r-1}}$$
which means that $t,s_1,\ldots,s_{r-1}$ satisfy  the same {\em braid relations}
as the Coxeter group $W(B_r)$ (actually the braid group of $G(d,1,r)$ is
isomorphic to that of $W(B_r)$, see \cite{BMR}), and that one gets a 
presentation of $G(d,1,r)$ by adding the {\em order relations} $t^d=1$ and
$s_i^2=1$.

For $G(de,e,r)$, to get a presentation of the braid group one needs to replace
the relation $ts_1ts_1=s_1ts_1t$ implied by the above diagram by the following
braid relations involving $t'$ and $s'_1$:
\begin{equation}\label{braidGde,e,r}
\begin{aligned}
t's'_1s_1&=s'_1s_1t',\\
 s'_1s_1s_2s'_1s_1s_2&=s_2s'_1s_1s_2s'_1s_1,\\
\underbrace{s_1t's'_1s_1s'_1s_1\ldots}_{e+1}&=
\underbrace{t's'_1s_1s'_1s_1\ldots}_{e+1}\\
\end{aligned}
\end{equation}
and one gets a presentation of $G(de,e,r)$ by adding the order relations
$t'^d=1$ and $s^{\prime 2}_1=s_i^2=1$.

\subsection{\label{RepGd1r}The representations of $G(d,1,r)$}

We   recall  from  \cite{ariki-koike}  explicit  models  over $K=\BQ(\zeta_d)$
for  the  irreducible
representations  of $G(d,1,r)$.  Let $\CL$  be the  set of  $d$-tuples $\bla =
(\lambda_0, \dots,\lambda_{d-1})$ of partitions with total size $r$. We denote
by  $\CT(\bla)$ the set of standard tableaux  of shape $\bla$, that is the set
of  tuples $\bT =  (T_0,\dots, T_{d-1})$ such  that $T_i$ is  a filling of the
Young diagram of $\lambda_i$ with numbers in $[1,r]$, with the conditions that
each  of these numbers appears  (exactly once) in one  of the filled diagrams,
and that they are increasing across the rows and columns of each $T_i$. For $m
\in [1,r]$ we let $\bT(m) = i$ if $m$ is placed in $T_i$.

To  $\bla \in  \CL$ we  associate the  $\BQ$-vector space $V^0(\bla)$ of basis
$\CT(\bla)$,   and  define  $V(\bla)  =  V^0(\bla)  \otimes  \BQ(\zeta_{d})$.
Explicit  formulas in \cite{ariki-koike} describe a representation $\rho_\bla$
of  $G(d,1,r)$ over $V(\bla)$, and the  $\rho_\bla$ provide a complete set of
representatives  for  the  irreducible  representations  of  $G(d,1,r)$  (see
\cite{ariki-koike}  corollary  3.14).  

$\bullet$ The matrix $\rho_\bla(t)$ is given by 
$\rho_\bla(t)\bT = \zeta_{d}^{\bT(1)} \bT$.

$\bullet$ The matrices $\rho_\bla(s_i)$ for $1 \leq i \leq n-1$ are rational in
the basis $\CT(\bla)$, thus belong to $\GL(V^0(\bla))$. To describe them
we introduce some notation: given a  tuple of tableaux $\bT$, let
$\bT_{i\leftrightarrow  i+1}$ be the tuple  obtained by exchanging the numbers
$i$  and $i+1$ in $\bT$  when this is still  a tuple of standard tableaux, and
$0$  otherwise (this last case will only  occur when the numbers $i$ and $i+1$
occur  in the same  tableau and on  either the same  line or the same column).
Then,  if $i$ and $i+1$ occur in the same tableau, define their axial distance
$a(i,i+1)$  to be  the distance  between the  diagonals where they occur (more
precisely,  if $i$  occurs at  coordinates $i_0,i_1$  and $i+1$ at coordinates
$j_0,j_1$  we  set  $a(i,i+1)=(i_0-i_1)-(j_0-j_1)$).  It  is  clear that for a
standard  tableau we always have $a(i,i+1)\ne 0$. Finally, take the convention
that   when  $i$  and  $i+1$   do  not  occur  in   the  same  tableau,  then
$a(i,i+1)=\infty$ so that $1/a(i,i+1)=0$. Then the formula is:
\begin{equation}\label{rhosi}
\rho_\bla(s_i)\bT= \frac 1{a(i,i+1)}\bT+(1+\frac 1{a(i,i+1)})
\bT_{i\leftrightarrow i+1}
\end{equation}

\subsection{\label{impr}The representations of $G(de,e,r)$}
As in \cite{ariki}, we will get the irreducible representations of $G(de,e,r)$
by Clifford theory, as components of the restriction of 
$\rho_\bla\in\Irr(G(de,1,r))$ to  $G(de,e,r)$.

We recall  that $\{t^e,s_1^t,s_1,\ldots,s_{r-1}\}$ are the
generating reflections for $G(de,e,r)$, and that the quotient
$G(de,1,r)/G(de,e,r)$ is cyclic, generated by $t$.
Let  $\chi$ be the generator of the group of linear characters of the quotient
which  maps $t$  to $\zeta_e$.  From the  formulae for  the action  of $t$ and
$s_i$,  we see that there is a permutation $\sigma(\bla)$ of the tuple $\bla$
such  that  $\chi\otimes\rho_\bla\simeq  \rho_{\sigma(\bla)}$;  it is given by
$\sigma(\bla)=(\lambda_d,\lambda_{d+1},\ldots,\lambda_{d+de-1})$   (where  the
indices  are taken $\pmod{de}$). Let 
$\langle  \sigma\rangle$  be the cyclic group generated  by $\sigma$, and let
$\centl=\langle\sigma^b\rangle$ be the subgroup which
stabilizes  $\bla$. By Clifford theory,
the  representation  $\rho_\bla$  restricts  to  $G(de,e,r)$  as  the  sum  of
$\#\centl=e/b$ distinct irreducible representations.

We  can actually define an operator  $S$ on $V(\bla)$ associated to $\sigma^b$
whose    eigenspaces    will    be    the    irreducible   constituents.   For
$\bT=(T_0,\ldots,T_{de-1})$     we     set    $S(\bT)=(T_{bd},T_{bd+1},\ldots,
T_{bd+de-1})$.  Since $\sigma^b(\bla)=\bla$, this is another tuple of tableaux
of shape $\bla$. It is easy to check that $S$ commutes to $\rho_\bla(s_i)$ and
to $\rho_\bla(s_1^t)$, and that $S\circ\rho_\bla(t)=\zeta_e^b\rho_\bla(t)\circ
S$.  It  follows  that  $S$  commutes  with  the action of $G(de,e,r)$ and its
eigenspaces  $V(\bla,\omega)  =  \Ker(S-  \omega)$  for  $\omega  \in  \mu_{\#
\centl}$  afford  irreducible  representations  of  $G(de,e,r)$.  We denote by
$\rho_{\bla,\omega}$  the representation  afforded by  $V(\bla,\omega)$; it is
clear       that       $$p_\omega:=\frac       1{\#\centl}\sum_{i=0}^{\#\centl
-1}\omega^{-i}S^i$$     is    the     $G(de,e,r)$-invariant    projector    on
$V(\bla,\omega)$.

From this we get, if we denote $\chi_{\bla,\omega}$ the character of
$\rho_{\bla,\omega}$:
$$\chi_{\bla,\omega}(g)=\frac 1{\#\centl}\sum_{i=0}^{\#\centl
-1}\omega^{-i}\Trace(\rho_\bla(g)S^i).$$

We  will now want to describe a model of $\rho_{\bla,\omega}$. We first recall
Ariki's  model (\cf. \cite[section  2]{ariki}), and then  introduce another
one, simpler for our purpose. Ariki chooses a basis of  $V(\bla,\omega)$ given
by the $p_\omega(\bT)$,
where he chooses for $\bT$ representatives of the $S$-orbits on $\CT$ given by
the  subset $\CT_0$ of  tuples of tableaux  which satisfy $\bT(1)<bd$. Setting
$\bT\om=p_\omega(\bT)$,  we then get formulas for the action of the generators
of $G(de,e,r)$ on the basis $\{\bT\om\}_{\bT\in\CT_0}$. We will write $\theta$
for $\zeta_e^b$ to simplify notations.
Using    that
$\rho_\bla(t)p_\omega=p_{\omega\theta}\rho_\bla(t)$, we get $\rho_{\bla}(t)
\bT\om=\zeta_{de}^{\bT(1)} \bT^{(\omega\theta)}$ from which we get
$$
\rho_{\bla,\omega}(t^e) \bT\om=\zeta_d^{\bT(1)} \bT\om.
$$
To  write  the  formula  for  $s_i$,  first  note  that  if $\bT\in\CT_0$ then
$\bT_{i\leftrightarrow  i+1}$ is not in general in $\CT_0$ when $i=1$. We have
to take $\bT'_{1\leftrightarrow 2}= S^{\lfloor\frac{\bT(2)}{bd}\rfloor}
\bT_{1\leftrightarrow 2}$ to get an element of $\CT_0$. Using that $p_\omega(
\bT'_{1\leftrightarrow 2})=\omega^{\lfloor\frac{\bT(2)}{bd}\rfloor}
p_\omega(\bT_{1\leftrightarrow 2})$, we get for $i>1$
$$
\rho_{\bla,\omega}(s_i)\bT\om= \frac 1{a(i,i+1)}\bT\om+(1+\frac 1{a(i,i+1)})
\bT\om_{i\leftrightarrow i+1}.
$$
and
$$
\rho_{\bla,\omega}(s_1)\bT\om= \frac 1{a(1,2)}\bT\om+(1+\frac 1{a(1,2)})
\omega^{-\lfloor\frac{\bT(2)}{bd}\rfloor}\bT\pom_{1\leftrightarrow 2}.
$$
Finally, using the above formulas we get
$$
\rho_{\bla,\omega}(s_1^t)\bT\om= \zeta_{de}^{\bT(1)-\bT(2)}
\rho_{\bla,\omega}(s_1)\bT\om.
$$

We now introduce a model which does not depend on a choice of representatives
$\CT_0$, using the equality $p_{\omega\theta\inv}\circ\rho_\bla(g)=
p_\omega\rho_\bla(\Ad t(g))$ which implies
$\chi_{\lambda,\theta^i}(g)=\chi_{\lambda,1}(\Ad(t^{-i})g)$.
We choose the same model as Ariki of $\rho_{\lambda,1}$, except that,
as the image by $p_1$ of a tuple $\bT$ is the same as that of all elements
of the same $S$-orbit, we take as basis elements the averages of $S$-orbits. 
We denote by $(\bT)$ the average of the $S$-orbit of $\bT$. We thus get a basis 
$\{(\bT)\}_{(\bT)\in \CT(\bla)/S}$ of $V(\bla,1)$. We take the same basis
for $V(\bla,\theta^i)$ and define
$\rho_{\bla,\theta^i}(g):=\rho_{\bla,1}(\Ad(t)^{-i} g)$.
We now have the following formulas for $\rho_{\bla,1}$:
\begin{equation}\label{rhobla1}
\begin{aligned}
\rho_{\bla,1}(t^e) (\bT)&=\zeta_d^{\bT(1)} (\bT),\\
\rho_{\bla,1}(s_i)(\bT)&=\frac 1{a(i,i+1)}(\bT)+(1+\frac 1{a(i,i+1)})
(\bT_{i\leftrightarrow i+1}) \text{ for any $i$}\\
\rho_{\bla,1}(s_1^{t^k})(\bT)&=\zeta_{de}^{k(\bT(1)-\bT(2))}
\rho_{\bla,1}(s_1)(\bT).
\end{aligned}
\end{equation}
(note that the above formulas make sense since both $\zeta_d^{\bT(1)}$ and
$\zeta_{de}^{\bT(1)-\bT(2)}$ are constant on a given $S$-orbit).

\section{The group $\Nbar$}
For $G$ irreducible, the  subgroup $\Nbar$ of $\Out(G)$ induced  by $N_{\GL(V)}(G)$ has been determined
in  \cite[3.13]{BMM}. For  the convenience  of the  reader we  recall here the
result.  When $K$ is a real field, $G$ is a finite Coxeter groups and elements
of  $\Nbar$ correspond to automorphisms of the Coxeter diagram. Such automorphisms
are traditionally denoted by an exponent, for example $\lexp 3 D_4$ denotes an
automorphism  of  order  $3$  of  the  group  $W(D_4)=G(2,2,4)$. We will use a
similar  notation for the other  cases. With one  exception $\lexp 2F_4$, all such
automorphisms  result from a  normal embedding of  $G$ into another reflection
group  of same  rank (the  case of  $\lexp 2F_4$  results from the embedding $
G_{28}=W(F_4)\subset  G_{31}$, which is  not normal). We  indicate at the same
time the relevant embedding.
$$
\begin{array}{|l|c|}
\hline 
\text{Automorphism}&\text{Embedding}\\
\hline 
\lexp iG(de,e,r) \text{where $i>1$, $i|de$}& G(de,e,r)\triangleleft G(de,1,r)\\
\lexp 3G(4,2,2)&G(4,2,2)\triangleleft G_6\\
\lexp 4G(3,3,3)&G(3,3,3)\triangleleft G_{26}\\
\lexp 3D_4& G(2,2,4)=W(D_4)\triangleleft G_{28}=W(F_4)\\
\lexp 2G_5&G_5\triangleleft G_{14}\\
\lexp 2G_7&G_7\triangleleft G_{15}\\
\lexp 2F_4&\\
\hline 
\end{array}
$$
This table needs a few comments on how it describes the group $\Nbar$.

The  automorphisms of $G(de,e,r)$  induced by $G(de,1,r)$  are induced by $\Ad
t$,   where  $t=\Diag(\zeta_{de},1,\ldots,1)$.  Let   us  determine  when  the
endomorphism  $\Ad t^i$ of $G=G(de,e,r)$ is inner. This happens if and only if
there    exists   $g\in   G$   such    that   $t^{-i}g$   is   scalar.   Since
$t^i=\Diag(\zeta_{de}^i,1,\ldots,1)$,  the  element  $g$  must  be of the form
$\Diag(\zeta_{de}^i  a,a,\ldots,a)$ for some  $a\in\mu_{de}$ and the condition
that    $g\in    G$    is    $a^r\zeta_{de}^i\in\mu_d$,   \ie.   $\zeta_{de}^i
\in\mu_d\cdot\mu_{de}^r=\mu_{de}^e\cdot\mu_{de}^r= \mu_{de}^{\gcd(e,r)}$, \ie.
we  get that $\Ad t^i$ is inner if and only if $i$ is multiple of $\gcd(e,r)$.
Thus the image of $\Ad t$ in $\Out(G)$ is of order $\gcd(e,r)$.

For $G=G(2,2,4)=D_4$ the group $\Nbar$ is generated by $\lexp3D_4$ and
$\lexp2G(2,2,4)$ and is isomorphic to $\Sgot_3$.

For $G=G(4,2,2)$ the group $\Nbar$ is generated by $\lexp3G(4,2,2)$ and
$\lexp2G(4,2,2)$ and is isomorphic to $\Sgot_3$.

For $G=G(3,3,3)$ the group $\Nbar$ is generated by $\lexp3G(3,3,3)$ and
the image in $\Out(G)$ of $\lexp4G(3,3,3)$, which is of order 2,
and is isomorphic to the alternating group $\Agot_4$.

The determination of $\Nbar$ has the following consequence.

\begin{proposition} \label{propfidirr} The natural action of
$\Out(G)$ on $\Irr(G)$ is faithful.
\end{proposition}
\begin{proof}
We note that an element of $\Out(G)$ which acts trivially on
$\Irr(G)$ belongs to $\Nbar$, as it fixes $\chi_V$. Assume
that $G = G(de,e,r)$, let $v = \gcd(e,r)$ and consider the $de$-tuple
of partitions $\bla = (\lambda_0,\dots,\lambda_{de})$
where $\lambda_i=\{r/v\}$ if $i$ is a multiple of $de/v$ and
$\lambda_i = \emptyset$ otherwise. With the notations of subsection
\ref{impr}, the tuple $\bla$ is stabilized by $\sigma^{e/v}$ which
corresponds to $(\Ad t)^{e/v}$, thus the restriction to $G$ of
$\rho_{\bla} \in \Irr G(de,1,r)$ has an inertia group of cardinality
$v$, hence the subgroup of $\Nbar$ induced by $\Ad t$
acts faithfully on its set of irreducible components.
It thus remains to consider the special cases mentioned above.
Inspection of the character tables, verifying that they are not left
invariant by the automorphisms in these special cases, concludes the proof.
\end{proof}

\section{\label{structout}The structure of $\Aut(G)$}
In this section we show that for an irreducible $G$, every automorphism is the
product  of a  central automorphism  and an  automorphism which  preserves the
reflections (except for $G=G(1,1,6)=\Sgot_6$). Since, for $G$ of rank at least
3, a non-trivial central automorphism does not preserve the reflections, this will prove
theorem \ref{auts}.

If  $G\subset\GL(V)$ be an  irreducible complex reflection  group, we may
identify $ZG$ to a subgroup of $\BC^\times$, thus $C$ is
formed     of     the     automorphisms     $\alpha_\chi$    of    the    form
$\alpha_\chi(g)=g\chi(g)\inv$  where  $\chi\in\Hom(G,\BC^\times)$  is a linear
character.  This defines in  general an endomorphism  of $G$ if  the values of
$\chi$ lie in $ZG$. it is an automorphism (thus a central automorphism) if and
only if $z\ne\chi(z)$ for all $z\in ZG-\{1\}$.

\subsection{The infinite series}
\begin{proposition} \label{propAutGdeer}
Any automorphism of a group $G(de,e,r)$ is  composed of an
automorphism which preserves the reflections with a central automorphism,
except for the case of $G(1,1,6)=\Sgot_6$ and for the case
of the non-irreducible group $G(2,2,2)\simeq\Sgot_2\times\Sgot_2$.
\end{proposition}

We first need a lemma.
\begin{lemma}
The group $D=D(de,e,r)$ of diagonal matrices in $G=G(de,e,r)$
is a characteristic subgroup of $G$ if and only if $G$ is not one of the
groups $G(2,2,2)$, $G(2,1,2)$, $G(4,2,2)$, $G(3,3,3)$, $G(2,2,4)$.
\end{lemma}
\begin{proof}
We  show first that  for $r\ge 5$  the group $D$  is the unique maximal normal
abelian subgroup  of  $G$.  Let  $E$  be  another  abelian  normal subgroup. Using the
semi-direct product decomposition $G=D\rtimes\Sgot_r$, the image $E/(E\cap D)$
of $E$ in $\Sgot_r$ is an abelian normal subgroup of $\Sgot_r$. Since the only
such subgroup of $\Sgot_r$ is trivial if $r\ge 5$, we get that $E/(E\cap D)=1$
meaning $E\subseteq D$ q.e.d.

We   now  assume  $r\le   4$.  If  $d=e=1$   then  $G=\Sgot_r$  and  $D=1$  is
characteristic. We thus assume now $de>1$. We will show that $D$ is the unique
normal  abelian subgroup  of maximal  order, excepted  for a few special cases
that  we  handle  separately.  As  before  we  exclude the non-irreducible and
abelian  case $G = G(2,2,2)$, for  which $\Out(G) = \GL_2(\mathbbm{F}_2)$ does
not  preserve $D\simeq \BZ/2$.  We note for  future use that  the order of the
center   $ZG$  is  the  gcd  of  the  reflection  degrees  of  $G$,  equal  to
$d\cdot\gcd(e,r)$,  and  that  $D$  is  of  order  $d^r  e^{r-1}$ (one chooses
arbitrarily   $r-1$  eigenvalues  in  $\mu_{de}$,   and  then  there  are  $d$
possibilities for the last one such that the product is in $\mu_d$).

For  $r=1$ we have $G=D$ thus there is nothing to prove. Let $E$ be an abelian
normal  subgroup of $G$  whose image in  $\Sgot_r$ is non-trivial.  For $1 < r
\leq  4$, the only non-trivial abelian normal subgroup $K$ of $\Sgot_r$ is transitive
and has order $r$. It follows that the diagonal matrices $D'$ in $E$ must have
equal eigenvalues, that is be in $ZG$. Thus the maximal such $E$ is the direct
product $ZG \times K$, and $\#E = r d\cdot\gcd(e,r)$. It is straightforward to
check  that  $\#E<\#D  =  d^r  e^{r-1}$  except  for  $G$  one  of $G(4,4,2)$,
$G(2,1,2)$,  $G(4,2,2)$, $G(3,3,3)$, $G(2,2,4)$. A  check shows that the outer
automorphism  of $G(2,1,2)$ does not preserves the diagonal matrices, but that
it   does   preserve   them   in   the   incarnation   $G(4,4,2)$;   and  that
$G(4,2,2),G(3,3,3),G(2,2,4)$ are genuine exceptions. 
\end{proof}

We are now in position to prove proposition \ref{propAutGdeer}, except for
the groups $G(2,1,2)$, $G(4,2,2)$, $G(3,3,3)$, $G(2,2,4)$, which will be
dealt with in the next section. We thus assume that $G$ is not
one of these exceptions, and in addition first assume $r \neq 6$.

Let $t$ and $s_1,\ldots,s_{r-1}$ be reflections as in section \ref{Gde,e,r}, such that
$G$  is generated by  $t^e, s_1^t,s_1,\ldots s_{r-1}$  (if $d=1$ we drop $t^e$
and  if $e=1$ we drop $s_1^t$). Let $\phi$ be an automorphism of $G$. Since we
are  not in  one of  the exceptions  of the  lemma, $\phi$  preserves $D$ thus
induces  an automorphism $\overline\phi$ of  $\Sgot_r$. Since we assumed $r\ne
6$,  the automorphism  $\overline\phi$ is  inner, of  the form $\Ad\sigma$ for
some  $\sigma\in\Sgot_r$, and lifting $\sigma$ to a permutation matrix in $G$,
we  see that up to an inner automorphism we may assume that $\phi$ induces the
identity  on $\Sgot_r$, that is,  preserves the shape of  a monomial matrix in
$G$.

We   will  denote  by  $M_\sigma$  the  permutation  matrix  corresponding  to
$\sigma\in\Sgot_r$. Thus $\phi(s_1)$     is      of     the     form
$\Diag(x_1,\ldots,x_r)M_{(1,2)}$;  the fact that  $\phi(s_1)$ is an involution
implies  that $x_1=x_2\inv$ and that $x_3,\ldots,x_n$ are signs. The fact that
$\phi(s_1)$   commutes  with   $\phi(s_3),\ldots,\phi(s_{r-1})$  implies  that
$x_3=x_4=\ldots=x_r$.      Thus     $\phi(s_1)$     is     of     the     form
$\Diag(\alpha_1\inv,\alpha_1,\epsilon,\ldots,\epsilon) M_{(1,2)}$ for  some
sign $\epsilon$ and some $\alpha_1\in\mu_{de}$.

Similarly  $\phi(s_i)$ is of the  form $\Diag(\epsilon, \ldots, \alpha_i\inv,
\alpha_i,\ldots,\epsilon)M_{(i,i+1)}$ where $\alpha_i$ is in $i$th position 
(the sign 
$\epsilon$ is  the same since $s_1$ and $s_i$ are conjugate). And $\phi(s_1^t)$
is of the form  $\Diag(\alpha^{\prime-1},  \alpha',\epsilon',\ldots,\epsilon')
M_{(1,2)}$.
We have actually $\epsilon=\epsilon'=1$  since  they occur only when $r>2$ in
which case $s_1^t  = (s_1)^{s_2  s_1 s_1^t  s_2}$ is conjugate to $s_1$.

Let $s=\Diag(1,\epsilon\alpha_1\inv,\epsilon^2(\alpha_1\alpha_2)\inv,\ldots,
\epsilon^r(\alpha_1\ldots\alpha_r) \inv)\in  G(de,1,r)$.  Then  $\Ad  s$
induces  an  automorphism of $G$ which
preserves   the   reflections,   and   the  composed  automorphism  $\phi'=\Ad
s\circ\phi$   satisfies   $\phi'(s_1)=\epsilon  s_1,\ldots,\phi'(s_r)=\epsilon
s_r$.      The     element      $\phi'(s_1^t)$     is      of     the     form
$\Diag(\alpha^{\prime\prime    -1},\alpha'',\epsilon,\ldots,\epsilon)$
for   some $\alpha''\in\mu_{de}$.

The  fact  that  $s_1s_1^t$  is  of  order  $de$  implies that $\alpha''$ is a
primitive  $de$-th  root  of  unity.  Thus  there exists an element $\gamma\in
\Gal(\BQ(\zeta_{de})/\BQ)$ such that  $\gamma(\alpha'')=\zeta_{de}$. 
Applying $\gamma$ to the matrices of $G$ gives an automorphism which preserves
the reflections,  and composing $\phi'$  with  that automorphism,  we  get  an
automorphism    $\phi''$    such    that    $\phi''(s_i)=\epsilon   s_i$   and
$\phi''(s_1^t)=\epsilon s_1^t$.

If $d=1$ we are finished: $\phi''$ is the central automorphism given by
$s\mapsto\epsilon s$ for any reflection.

Assume  now $d>1$. Then $\phi''(t^e)$ is a diagonal matrix which commutes with
$s_2,\ldots,s_{r-1}$  thus is of the  form $z\Diag(\zeta,1,\ldots,1)$ for some
$z\in  ZG$ and some $\zeta\in\mu_{de}$. Since it  is of order $d$, $z$ must be
of order $d$ and $\zeta$ be a primitive $d$-th root of unity.

Assume first that $e=1$. Then we do not need $s_1^t$ among the generators, and
composing $\phi''$ with the automorphism induced by some appropriate element of
$\Gal(\BQ(\zeta_{de})/\BQ)$ we may get a $\phi'''$ such that
$\phi'''(s_i)=\epsilon  s_i$ and $\phi'''(t)=z'  t$ for some  $z'\in ZG$, thus
$\phi'''$ is a central automorphism.

Finally  assume that $e$ and  $d$ are not $1$.  Then we use the third relation
$\underbrace{s_1t^es_1^ts_1s_1^ts_1\ldots}_{e+1}=
\underbrace{t^es_1^ts_1s_1^ts_1\ldots}_{e+1}$   in   \ref{braidGde,e,r}  which
implies  that  $\zeta=  \zeta_d$,  thus  $\phi$  is  the  central  automorphism
$s_i\mapsto\epsilon s_i$ and $t^e\mapsto z t^e$.

The  same  method  will  apply  when  $r=6$  and $de>1$ if we can exclude
the possibility for $\overline\phi$ to be a non-inner automorphism of
$\Sgot_6$.  Up  to  conjugation  by  a  permutation matrix, we may assume that
$\overline\phi$  sends  $(1,2)$  to  $(1,2)(3,4)(5,6)$.  Let  
$p=M_{(1,2)(3,4)(5,6)}$. We  have
$\phi(C_G(s_1)\cap  D)  =  C_G(\phi(s_1))\cap  D  =  C_G(p) \cap D$, the first
equality  since $D$ is a characteristic subgroup.  and the second since $G = D
\rtimes  \Sgot_6$ and $D$  is abelian. It  follows that $C_G(s_1)  \cap D$ and
$C_G(p)  \cap D$ have  the same cardinality.  One easily gets $\#(C_G(s_1)
\cap  D) = (de)^4 d$. On the other hand, a matrix in $C_G(p) \cap D$ is
uniquely  determined by (diagonal) coefficients  $(a,b,c) \in \mu_{de}^3$ such
that $c^2 \in a^{-2} b^{-2} \mu_d$. Since an element in $\mu_{de}$ has at most
two  square  roots,  it  follows  that  $\#(C_G(p) \cap D )\leq 2d
(de)^2$. But $de \geq 2$ implies $ 2d(de)^2 < (de)^4 d$, a contradiction.
\subsection{The other cases}
We  now  handle  the  exceptional  groups  as  well as ones left over from the
previous   subsection.  We  first  establish   a  reduction  result.  We  call
hyperplanes of $G$ the hyperplanes of $V$ which are invariant subspaces of the
reflections  in  $G$.  Up  to  a scalar  there  is a unique $G$-invariant scalar
product  on $V$, to which the  orthogonality statement in the next proposition
refers.

\begin{proposition}
Let  $\phi$ be an automorphism of an  irreducible reflection group $G$
which maps $R$ inside $R\cdot ZG$, where $R$ is the set of reflections of
$G$.  Then if  either $G$ is of rank $\ge  3$ or  $G$ has  the property  that
two orthogonal hyperplanes  are  conjugate then  $\phi\in C\cdot A$.
\end{proposition}
\begin{proof}
The  group $G$ admits a presentation by a set $S$ of distinguished reflections
(\ie,  reflections  $s$  with  eigenvalue  $e^{2i\pi/e_s}$, where $e_s$ is the
order  of the centralizer  in $G$ of  the reflecting hyperplane  of $s$), with
relations  the order  relations $s^{e_s}=1$  and braid  relations $w=w'$ where
$w,w'\in  S^*$ are words of the same  length; here the braid relations present
the  braid group of  $G$ (see \cite[theorem  0.1]{bessis2}). We note that two
distinguished  reflections  are  conjugate  if  and  only  if their reflecting
hyperplanes  are conjugate. The  words $w$ and  $w'$ have the  same number of
elements  in each class of reflections,  since the braid relations present the
braid  group, whose  abelianized is  the free  abelian group  on the conjugacy
classes of hyperplanes.

For  any $s\in S$ we define a reflection $f(s)$ and an element $z(s)\in ZG$ by
$\phi(s)=f(s)z(s)$.  When  $\dim  V\ge  3$,  $f(s)$  and  $z(s)$  are uniquely
determined,  but when $\dim V=2$  there might be two  choices. We then make an
arbitrary  such choice,  but we  ask that  if $s,s'\in  S$ are conjugate, \ie.
there  exists $g\in  G$ such  that $s'=\lexp  gs$ then  $z(s')=z(s)$, which is
possible  since $\phi(s')=\lexp{\phi(g)}f(s) z(s)$ and $\lexp{\phi(g)}f(s)$ is
a  reflection. We now  extend $f$ and  $z$ to the  free monoid $S^*$; we claim
that  this  gives  well-defined  group  homomorphisms  from $G$. Indeed, since
(obviously)  for any  $w\in S^*$  we have  $\phi(w)=f(w)z(w)$, it is enough to
show  that $z$  induces a  linear character  of $G$.  From the  fact that each
conjugacy of reflections occurs as many times in $w$ and $w'$ it is clear that
$z(w)=z(w')$.    And   from    $\phi(s^{e_s})=1=f(s)^{e_s}z(s)^{e_s}$,   since
$f(s)^{e_s}$   has   at   least   one   eigenvalue   $1$,   we  conclude  that
$f(s)^{e_s}=z(s)^{e_s}=1$.

If  $\dim V\ge 3$, then $f$ is an  automorphism: to show that, it is enough to
show that $f$ is surjective on $R$, and for that it is enough to show that $f$
is    injective   on   $R$.    If   $s,s'\in   R$    and   $f(s)=f(s')$   then
$\phi(ss^{\prime-1})\in  ZG$ which implies $ss^{\prime-1}\in ZG$ which implies
$s=s'$  since $\dim V\ge 3$. Thus if $\dim  V\ge 3$, we get that $\phi$ is $f$
composed with a central automorphism (given by $s\mapsto s z(f\inv(s))$).

It  remains the  case $\dim  V=2$ and $G$ has the property that two orthogonal
hyperplanes are conjugate for  which we  give  a more complicated
argument  to prove that $f$ is injective. It is enough to show that some power
of  $f$ is injective.  Since $\phi$ and  $f$ both map  $ZG$ into itself, it is
easy  to show by induction on $r$  that for any $r$ the endomorphisms $\phi^r$
and $f^r$ still differ by a linear character with values in $ZG$. Since $\phi$
is  an automorphism,  it is  of finite  order; taking  for $r$ this order, and
replacing $f$ by $f^r$, we see that we may assume $\phi=\Id$.

We  thus have $f(g)=g z(g)$ for some  $ZG$-valued linear character $z$, and we
want  to see that $f$ is injective. We  will show that $f^2=\Id$. It is enough
to  show that for a reflection $r$. We have two cases to consider: if $f(r)=r$
then $f^2(r)=r$. Otherwise, we must have $z(r)= \zeta\inv\Id$ where $\zeta$ is
the non-trivial eigenvalue of $r$, and if $\langle,\rangle$ is a $G$-invariant
scalar  product on  $V$, $f(r)$  is a  reflection with  non-trivial eigenvalue
$\zeta\inv$  and  hyperplane  orthogonal  to  that  of  $r$.  Since  two  such
hyperplanes  are conjugate and $f(r)\inv$ has the same eigenvalue as $r$, they
are          conjugate.          Then          $z(f(r))=z(r)\inv$         thus
$f(f(r))=f(r)z(f(r))=rz(r)z(r)\inv=r$ and we are done.
\end{proof}

By  enumerating the hyperplanes in \CHEVIE (see \cite{chevie}), it  is easy to see that orthogonal
hyperplanes  are  conjugate  for  all  2-dimensional exceptional groups except
$G_5$ and $G_7$.

The  preceding proposition is sufficient to handle the exceptional
groups $G_i$ for  $i\ge 12$. Indeed, an easy computer
check shows that any conjugacy class $c$ such that there exists a conjugacy
class $c'$ satisfying
\begin{itemize}
\item $\#c' = \#c$
\item $c'$ contains a reflection $s'$
\item $c$ contains an element with the same order as $s'$
\end{itemize}
contains  the product  of a  reflection by  an element  of the center. 

Finally, for the groups  $G_4$ to $G_{12}$, as well  as for the exceptions 
in the infinite series $G(2,1,2), G(4,2,2),
G(3,3,3)$ and $G(2,2,4)$ we use the techniques of section \ref{exceptional} to
determine  all automorphisms  and to  check that  they are  the product  of an
element of $C$ and an automorphism which preserves the reflections.

\section{Reflection representations}

By a reflection representation of a reflection group $G$, we mean a faithful
representation such that the image of a reflection of $G$ is a reflection.
The corresponding characters are called reflection characters.
Here we prove theorem \ref{thref}, whose statement we recall. 

\begin{theorem}  \label{refl} Let  $G\subset\GL(V)$ be  an irreducible complex
reflection group, where $V$ is a vector space over the field of definition $K$
of  $G$. Then the reflection characters of $G$ are the transformed of $\chi_V$
by the action of $\Gal(K/\BQ)$.
\end{theorem}
Before going on with the proof, we remark that this proves indeed that
any automorphism which preserves the reflections is a Galois automorphism
attached to some $\gamma\in\Gal(K/\BQ)$, since it sends $\chi_V$ to some 
$\gamma(\chi_V)$.
\begin{proof}[of \ref{refl}]
The theorem is proved by a computer check for the exceptional reflection groups,
by looking at the character tables. Both the action of the Galois group, and
which representations are faithful and send reflections to reflections are
readily seen from the character tables.

We now look at the groups $G(de,e,r)$.
We first look at the case $e=1$.  The result is clear if $G$ is cyclic so
we assume $r>1$, and use the notations of section \ref{RepGd1r}.

We first describe the action of $\Gal(K/\BQ)$ on the representations
of $G(d,1,r)$.

\begin{lemma}\label{gal on repGd1r}
For $i\in\BZ$ prime to $d$, let $\gal i\in \Gal(\BQ(\zeta_d)/\BQ)$ be the 
element
defined  by  $\gal  i(\zeta_d)  =  \zeta_d^i$.  Then the partition tuple $\gal
i(\bla)$ such that $\gal i\circ\rho_\bla\simeq\rho_{\gal i(\bla)}$ is given by
$\gal i(\bla)=(\lambda_0,\lambda_i,\lambda_{2i},\ldots, \lambda_{(d-1)i})$
(where the indices are taken $\pmod d$).
\end{lemma}
\begin{proof}
This is immediate using the model given in \ref{RepGd1r}. The action of
$\gal{i}$ is trivial on $\rho_\bla(s_i)$ and maps $\rho_\bla(t)$ to
$\rho_\bla(t)^i$, from which the result follows.
\end{proof}

We now determine which $\rho_\bla$ are faithful reflection representations.
\begin{lemma}\label{refl for Gd1r}
The faithful reflection representations of $G(d,1,r)$ are
the  $\rho_{\bla^{(i)}}$ for $i$ prime to $d$, where 
$\lambda^{(i)}_0=\{r-1\}$ and  $\lambda^{(i)}_i=\{1\}$. 
\end{lemma}
\begin{proof}
Let us determine when $\rho_\bla(t)$ is a reflection. First, notice that there
is  only one standard tableau corresponding  to a given young  diagram and with a given
content  if and only  if the diagram  is a column  or a row.  For a partition
$\lambda$  we denote  by $|\lambda|$  the size  of $\lambda$  (the sum  of the
lengths  of  its  parts), and call height of $\lambda$ the number of
its parts, namely the number of rows in the corresponding diagram.
Since $\bT\in\CT(\bla)$ are eigenvectors of
$\rho_\bla(t)$  with eigenvalue
$\zeta_{d}^{\bT(1)}$,  there must be exactly  one $\bT\in\CT(\bla)$ such that
$\bT(1)\ne  0$. This means that  there must be exactly  one $i\ne 0$ such that
$|\lambda_i|\ne 0$, with a Young diagram a line or a column. We must also have
that  $\lambda_0$  is  a  line  or  a  column;  if  $|\lambda_0|=0$  then  the
representation  is of dimension 1 and cannot  be faithful since $r>1$. Also if
$|\lambda_i|>1$  then $r\ge  3$ and  there are  at least  two $\bT$  such that
$\bT(1)=i$. So we must have $\lambda_i=\{1\}$. Finally, note that $i$ must
be  prime to $d$ otherwise $\rho_\bla$  is not faithful (since $\rho_\bla(t)$
does not have order $d$). To summarize, we have $\lambda_0=\{r-1\}$ or
$\{1^{r-1}\}$ and $\lambda_i=\{1\}$ for some $i$ prime to $d$.

Let  us now analyze when $\rho_\bla(s_i)$ is a reflection (or equivalently all
$\rho_\bla(s_i)$   for  $i=1,\ldots,r-1$  are   reflections,  since  they  are
conjugate).  The formula \ref{rhosi} shows  that $\rho_\bla(s_i)$ in the basis
$\CT(\bla)$  is  block-diagonal,  with  diagonal  blocks  indexed by the basis
elements $\bT$ and $\bT_{i\leftrightarrow i+1}$, which are thus of size 1 when
$\bT_{i\leftrightarrow i+1}=0$ and of size 2 otherwise.

\begin{enumerate}
\item   If  $\bT_{i\leftrightarrow  i+1}\ne  0$  the  corresponding  block  is
$\begin{pmatrix}\frac      1{a(i,i+1)}&1+\frac      1{a(i,i+1)}\\      1-\frac
1{a(i,i+1)}&\frac{-1}{a(i,i+1)}  \end{pmatrix}$ whose eigenvalues  are $1$ and
$-1$.
\item  If $\bT_{i\leftrightarrow i+1}=0$,  then $i$ and  $i+1$ are in the same
tableau,  and in the same line or column. We have that $\bT$ is an eigenvector
of  $\rho_\bla(s_i)$, for  the eigenvalue  $\frac 1{a(i,i+1)}=1$  when $i$ and
$i+1$ are in the same line and $\frac1{a(i,i+1)}=-1$ when $i$ and $i+1$ are in
the same column. 
\end{enumerate}

For $\rho_\bla(s_i)$ to be a reflection there cannot be more than two $j$ with
$|\lambda_j|\ne 0$ otherwise case (i) would occur more than once (with $i$ and
$i+1$  in  different  tableaux).  Similarly,  if  there  are two $j$ such that
$|\lambda_j|\ne  0$  then  the  two  Young  diagrams  must be lines or columns
otherwise case (i) would occur more than once. Actually, the two diagrams must
be  lines  since  if  one  is  a  column  of  height  $>1$ then case (ii) with
eigenvalue $-1$ has to occur at least once. And for $\rho_\bla(s_{r-1})$ to be
a  reflection one of the lines has to  be of length $1$, otherwise by removing
the  last square  of each  line we  can still  fill the  rest in  at least two
different ways which leads to too many instances of (i).

We now look at the case when only one $|\lambda_j|\ne 0$. The diagram must not
be  a line or a column otherwise $\rho$  is not faithful. For $i=r-1$ case (i)
will  occur more than once if the diagram  has more than two corners. If there
are  two corners, when they are  removed the rest must be  a line or a column;
this  means the diagram must be a hook; and the height of the hook has to be 2
otherwise  case (ii) with eigenvalue $-1$ will  occur. Finally if there is one
corner  the  diagram  is  a  rectangle;  this  rectangle has to be $2\times 2$
otherwise  by removing  the corner  one would  get a  diagram with two corners
which  is  not  a  hook  which  would  lead  to  too many eigenvalues $-1$ for
$s_{r-2}$.

To  summarize,  we  must  have  $i\ne  j$  such  that  $\lambda_i=\{r-1\}$ and
$\lambda_j=\{1\}$, or only one $i$ with $|\lambda_i|\ne 0$ with diagram either
a hook of height 2 or a $2\times 2$ square (the last for $r=4$).

So we are left with the cases in the statement of the theorem for the images
of all reflections to be reflections, and the representation to be faithful.
\end{proof}

The  reflection  representation  of  $G(d,1,r)$  is $\rho_{\bla^{(1)}}$ and by
lemma \ref{gal on repGd1r} we have $\gal
i\circ\rho_{\bla^{(1)}}=\rho_{\bla^{(i)}}$  which  together  with  lemma
\ref{refl for Gd1r} proves the theorem for $G(d,1,r)$.

Let us look now at $G(de,e,r)$. 

\begin{lemma}
The  faithful irreducible  reflection representations  of $G(de,e,r)$  are the
restriction   to   $G(de,e,r)$   of   the   faithful   irreducible   reflection
representations of $G(de,1,r)$.
\end{lemma}
\begin{proof}
We first notice that a representation which is
not  the restriction  of an  irreducible representation  of $G(de,1,r)$ is not
faithful.  To prove this, we use that  the central character of $\rho_\bla$ is
given   by   $\omega_{\rho_\bla}(z)=\zeta_{de}^{\sum_{i=0}^{d-1}i|\lambda_i|}$
(this  can be deduced from  the construction of $\rho_\bla$  as induced from a
generalized  Young subgroup --- see  \eg. \cite[p. 93--106]{Zelevinsky}). Here
$z$  is the generator of the center of $G(de,1,r)$, which can be identified to
the scalar $\zeta_{de}\in\GL(\BC^r)$. If $\rho_\bla$ restricted to $G(de,e,r)$
is  not irreducible, there  is a divisor  $b$ of $e$  such that for all $i$ we
have $\lambda_{i+de/b}=\lambda_i$ (then each $\lambda_i$ occurs $b$ times thus
$b$  divides $\gcd(e,r)$).  It follows  that $\sum_{i=0}^{d-1}i|\lambda_i|$ is
divisible  by $b$ thus  $\rho_\bla$ is not  faithful on the  center. It is not
even  faithful on the center of $G(de,e,r)$  since that center is generated by
$z^{\frac e{\gcd(e,r)}}$ and $b$ still divides the order of that element.

Thus the representation we are looking at is the restriction to $G(de,e,r)$ of
some   $\rho_\bla$.  We  use  the  analysis  of  when  $\rho_\bla(s_i)$  is  a
reflection.  Notice first that  the representations with  only one diagram are
not faithful, since by formula \ref{rhobla1} $s_1$ and $s_1^t$ have same image
since  for any $\bT$ we have $\bT(1)=\bT(2)$. Thus we are in the case where we
have  $i\ne j$ such that  $\lambda_i=\{n-1\}$ and $\lambda_j=\{1\}$. Actually,
for     the    representation     to    be     faithful    we     must    have
$\rho_\bla(s_1^{t^k})\ne\rho_\bla(s_1)$  for  any  $k$  prime  to  $de$,  thus
$\bT(2)-\bT(1)$ must be prime to $de$, \ie. $j-i$ must be prime to $de$.

Analyzing  when $\rho_\bla(t^e)$ is  a reflection along  the same lines we did
for  $\rho_\bla(t)$ we  find that  $j\equiv 0\pmod  d$. Since  $\rho_\bla$ and
$\rho_{\sigma(\bla)}$  have same  restriction to  $G(e,e,r)$ where $\sigma$ is
the  ``shift'' by $d$ as  in section \ref{impr}, we  may assume $j=0$.
Thus we have indeed the same representations as in
\ref{refl for Gd1r}, which proves the lemma.
\end{proof}

The theorem for $G(de,e,r)$ follows immediately from the above lemma since
the reflection representation of $G(de,e,r)$ is the restriction of the
reflection representation of $G(de,1,r)$, and the Galois action commutes with
the restriction.
\end{proof}

\section{\label{exceptional} Theorems \ref{A} and \ref{xx} for the 
exceptional groups}
To  verify \ref{xx}  for the  exceptional reflection  group, we  have used the
\GAP\  package \CHEVIE\  (see \cite{chevie}),  which contains information about
them, including:
\begin{itemize}
\item Their realization via a permutation representation (on a set of
``roots'' in $V$).
\item Their character table.
\end{itemize}

The  character  tables  allow  us  to  determine  the permutation of $\Irr(G)$
effected  by an  element of  $\Gal(K/\BQ)$. The  method we  have used to check
\ref{A}  is to compute the group of  outer automorphisms of $G$ which preserve
the  set of reflections,  and then to  find a subgroup  which induces the same
permutation group on $\Irr(G)$ as $\Gal(K/\BQ)$. We recall
from proposition \ref{propfidirr} that
an element of $\Out(G)$ is completely determined
by the permutation it induces on $\Irr(G)$.

To compute automorphisms, we use a presentation for $G$, by means of a diagram
(see  \cite{BMR} and \cite{bessis-michel});  if $S$ is  the generating set for
this   presentation,  to  enumerate  the   automorphisms  which  preserve  the
reflections,  we enumerate up to inner automorphisms all tuples of reflections
of  $G$ which satisfy the  relations for $G$, and  generate $G$.

If  we want to  enumerate all automorphisms,  we have to  extend the search to
other  conjugacy  classes  which  have  the  same  size  and same order as the
conjugacy  class  of  a  reflection.  We  only  need  to  compute the whole of
$\Aut(G)$  for $2$-dimensional groups, to meet  the requirements at the end of
section \ref{structout}, and for these groups the process above is quite fast.

To find directly the orbits of $G$ on the tuples of reflections is a too large
computation  for the larger groups, such  as $G_{34}$ which has about $4.10^7$
elements;  thus we use a  recursive process: we choose  the image of the first
generator  $s$ (which amounts  to choosing a  representative for each orbit of
$G$  on the  reflections); then  we find  the $C_G(s)$-orbits  of tuples which
satisfy  the relations for $S-\{s\}$ (by  a recursive process) and then select
those  that have  the required  relations with  $s$; and  finally weed out the
tuples which do not generate $G$. This allows us compute $\Abar$
in at most  a few hours even   for  the  largest  groups.

To find the  image of $\ibar_\chi(\Gamma)$ in $\Abar$ is  rather easy since
for the exceptional groups
the  morphism $\Gamma\xrightarrow\ibar\Abar$ is either surjective or has an
image
of  index 2 (the last occurs exactly for $G_5$, $G_7$ and $G_{28}=F_4$, which,
in  agreement with theorem \ref{A->Gal surjective}, are also  the only 
exceptional groups where $\Nbar\ne 1$; in these cases, with the notations of 
table \ref{titilde}, $\Nbar$
induces  respectively the  automorphism $(s,t)\mapsto(t,s)$,  the automorphism
$(s,t,u)\mapsto(s,u\inv,t\inv)$, the diagram automorphism of $F_4$).

To  check \ref{xx}, we must lift $\ibar$  to the group $\Aut(G)$, that is find
representatives  in $\Aut(G)$ of  elements of the  image  of  $\ibar$  which 
satisfy the relations for
$\Gamma$.  Since for exceptional groups $\Gamma$ is  an abelian $2$-group,  
these relations are order
relations, plus commutation relations.

Few  representatives of a given automorphism turn out to have the right order,
so  we first compute all representatives which  have the right order, and then
it is an easy job to pick among them elements which commute (this procedure of
course  fails for $G_{27}$  for which we  must construct $\itilde$ rather than
$\iota$;  in this case we  take $K'=\BQ(\zeta_{15})$; the group $\Gal(K'/\BQ)$
is  generated by $\cc$ and $\gal{7}$, this last element being of order $4$ (it
is of order 2 in $\Gal(K/\BQ)$) ; we can map these to automorphisms of order 2
and  4 respectively commuting to each other  --- but it  is impossible to map
$\gal{7}$ to an automorphism of order 2).

The  only computational problem in this procedure is to compute the product
of  two automorphisms;  since an  automorphism $\varphi$ is  represented by the 
image $\varphi(S)$,  to do this we need an expression as words in $S^*$
of the elements of $\varphi(S)$.
For  all  groups  but  $G_{34}$  this  was solved by enumerating all
elements  of $G$ by standard  permutation group algorithms using  a base and a
strong   generating  set.  For   $G_{34}$  we  managed   by  only  considering
automorphisms  which  extend  the  one  for  $G_{33}$; since the normalizer of
$G_{33}$   in  $G_{34}$  is  $G_{33}$  times  the  center  of  $G_{34}$,  such
automorphisms  map the last generator of  $G_{34}$ to a conjugate by $G_{33}$,
which makes the computation feasible.

Table \ref{titilde} summarizes our results for the exceptional reflection groups; it
gives  the diagram  for $G$,  and then  gives the  value of  $K$ and describes
$\iota$ in terms of the generators $S$ given by the diagram.

The nodes of the diagrams are labeled by the elements of $S$; we give to these
elements  names in the list $s, t, u, v,  w, x$ in that order. A number inside
the  node gives the order of the  corresponding reflection when it is not $2$.
To  express the braid  relations, we use  the conventions for Coxeter diagrams:
two  nodes  which  are  not  connected  correspond to commuting reflections. A
single  line  connecting  two  nodes  $s$  and $t$ corresponds to the relation
$sts=tst$, a double line to $stst=tsts$, a triple line to $ststst=tststs$, and
a   line  labeled  with  the  number  $n$  corresponds  to  a  braid  relation
$\underbrace{stst\cdots}_{n\text{  terms  }}= \underbrace{tsts\cdots}_{n\text{
terms }}$.

These  conventions are extended  by additional ones;  a circle joining 3 nodes
$s$,  $t$ and $u$, with the number $n$ inside (when there is no number inside,
3   is  to  be  understood)  corresponds  to  a  clockwise  circular  relation
$\underbrace{stust\cdots}_{n\text{ terms }}= \underbrace{tustu\cdots}_{n\text{
terms  }}=  \underbrace{ustus\cdots}_{n\text{  terms  }}$;  for  instance, the
relation  for $G_7$, $G_{11}$, $G_{19}$ and the generators $s,t,u$ of $G_{31}$
is  $stu=tus=ust$.  Similarly  the  $6$  inside  the triangle for $G_{33}$ and
$G_{34}$  means that in addition to length  $3$ braid relations implied by the
sides of the triangle, there is the circular relation $tuvtuv=uvtuvt=vtuvtu$.

The  double line in $G_{29}$ expresses a length $4$ braid relation between $v$
and  $ut$:  $vutvut=utvutv$.  Let us mention that the exceptional characters
of $G_{29}$ are the 6-dimensional rational ones which appear in the 10th symmetric
power $S^{10}(V)$.

The  triangle  symbol  for $G_{24}$ and $G_{27}$
corresponds  to the relation $s(uts)^2=(stu)^2t$. Finally, the groups $G_{13}$
and  $G_{15}$ have more complicated relations  which are spelled out below the
diagram.

The  groups $G_{28}=F_4$,  $G_{35}=E_6$, $G_{36}=E_7$  and $G_{37}=E_8$, which
are  rational, do not appear  in the table since  our theorems are trivial for
them.

The map $\iota$ is specified by giving for generators $\gamma$ of $\Gamma$ the
tuple  $\iota(\gamma)(S)$  (where  $S$  is  in  the  order $s,t,u,v,w,x$). The
elements  of $\Gamma$  are denoted  by $\cc$  for complex  conjugation and by
$\gal  i$  for  the  following  element:  let  $\BQ(\zeta_n)$  be the smallest
cyclotomic  field containing $K$; then  $\gal i$ is the  restriction to $K$ of
the   Galois  automorphism   of  $\BQ(\zeta_n)$   which  sends   $\zeta_n$  to
$\zeta_n^i$.

\font\hugecirc = lcircle10 scaled 1950
\font\largecirc = lcircle10 scaled  \magstep 3
\font\mediumcirc = lcircle10 scaled \magstep 2
\font\smallcirc = lcircle10 
\def\BBcirc{$\hbox{\hugecirc\char"6E}$}
\def\Bcirc{$\hbox{\largecirc\char"6E}$}
\def\sBcirc{$\hbox{\mediumcirc\char"6E}$}
\def\sbcirc{$\hbox{\smallcirc\char"6E}$}
\def\ucirc{$\hbox{\mediumcirc\char"13}$}
\def\lcirc{$\hbox{\mediumcirc\char"12}$}
\def\cnode#1{{\kern -0.4pt\bigcirc\kern-7pt{\scriptstyle#1}\kern 2.6pt}}
\def\node{{\kern -0.4pt\bigcirc\kern -1pt}}
\def\snode{{\kern -0.4pt{\scriptstyle\bigcirc}\kern -1pt}}
\def\tbar#1pt{{\rlap{\vrule width#1pt height1pt depth0pt}
           \rlap{\vrule width#1pt height3pt depth-2pt}
\vrule width#1pt height5pt depth-4pt}}
\def\overmark#1#2{\kern -1.5pt\mathop{#2}\limits^{#1}\kern -2pt}
\def\trianglerel#1#2#3{
 \nnode#1\bar14pt\kern-13pt\raise7.5pt\hbox{$\displaystyle \underleftarrow 6$}
  \kern 2pt\nnode#2
 \kern-29pt\raise9.5pt\hbox{$\diagup$}
 \kern -2pt
 \raise17.5pt\hbox{$\node$\rlap{\raise 2pt\hbox{$\kern 1pt\scriptstyle #3$}}}
 \kern -1pt\raise9.5pt\hbox{$\diagdown$}
 \kern 3pt
}
\def\vertbar#1#2{\rlap{$\nnode{#1}$}
                 \rlap{\kern4pt\vrule width1pt height17.3pt depth-7.3pt}
\raise19.4pt\hbox{$\node$\rlap{$\kern 1pt\scriptstyle#2$}}}
\def\stacktwo#1#2{\genfrac{}{}{0pt}{0}{#1}{#2}}
\def\stackthree#1#2#3{\genfrac{}{}{0pt}{0}{#1}{\genfrac{}{}{0pt}{0}{#2}{#3}}}
\begin{longtable}[c]{|>{$}c<{$}|>{$}c<{$}>{$}c<{$}>{$}c<{$}|}
\caption{The morphism $\itilde$}\label{titilde}\\
\hline \mbox{Group}&\mbox{Diagram}&\mbox{Field}&\mbox{$\itilde$}\\
\hline
\endfirsthead
\hline \mbox{Group}&\mbox{Diagram}&\mbox{Field}&\mbox{$\itilde$}\\
\hline
\endhead
\endfoot
\hline
\endlastfoot
\vrule width0pt height 0.5cm depth 0.5cm
G_4&\ncnode3s\bar20pt\ncnode3t&\BQ(\zeta_3)&\cc\mapsto(s\inv,\lexp st\inv)\\
\hline
\vrule width0pt height 0.5cm depth 0.5cm
G_5&\ncnode3s\dbar20pt\ncnode3t&\BQ(\zeta_3)&\cc\mapsto(s\inv,\lexp st\inv)\\
\hline
\vrule width0pt height 0.8cm depth 0.5cm
G_6&\nnode s\tbar20pt\ncnode3t&\BQ(\zeta_{12})&
\stacktwo{\cc\mapsto(s\inv,t\inv),}{\gal 7\mapsto(s^{t\inv st},t)}\\
\hline
\vrule width0pt height 0.8cm depth 0.7cm
G_7& {\scriptstyle s}\node\kern13.5pt\raise2.5pt\hbox{$\Bcirc$}
  \kern-15.5pt\rlap{\raise10pt\hbox{$\cnode3\kern 1pt\scriptstyle t$}}
  \lower10pt\hbox{$\cnode3\kern 1pt\scriptstyle u$}&\BQ(\zeta_{12})&
  \stacktwo{\cc\mapsto(s\inv,t\inv,(u\inv)^s),}
  {\gal 7\mapsto(s^{ut},t,u)}\\
\hline
\vrule width0pt height 0.5cm depth 0.5cm
G_8&\ncnode4s\bar20pt\ncnode4t&\BQ(\ii)&\cc\mapsto(s\inv,t\inv)\\
\hline
\vrule width0pt height 0.8cm depth 0.5cm
G_9&\nnode s\tbar20pt\ncnode4t&\BQ(\zeta_8)&
\stacktwo{\cc\mapsto(s\inv,t\inv),}{\gal 5\mapsto(s^{tts},t)}\\
\hline
\vrule width0pt height 0.8cm depth 0.5cm
G_{10}&\ncnode3s\dbar20pt\ncnode 4t&\BQ(\zeta_{12})&
\stacktwo{\cc\mapsto((s\inv)^{t\inv},t\inv)}{\gal 7\mapsto(s^{ttss},t\inv)}\\
\hline
\vrule width0pt height 0.85cm depth 0.85cm
G_{11}&
 {\scriptstyle s}\node\kern13.5pt\raise2.5pt\hbox{$\Bcirc$}
  \kern-15.5pt\rlap{\raise10pt\hbox{$\cnode3\kern 1pt\scriptstyle t$}}
  \lower10pt\hbox{$\cnode 4\kern 1pt\scriptstyle u$}
  &\BQ(\zeta_{24})&
  \stackthree{\cc\mapsto(s\inv,t\inv,(u\inv)^s),}
  {\gal {13}\mapsto(s^{tsut},t,u),}
  {\gal {19}\mapsto(s,t,(u\inv)^{tu\inv s})}\\
\hline
\vrule width0pt height 0.8cm depth 0.7cm
G_{12}&
 {\scriptstyle s}\node\kern13.5pt\raise2.5pt\hbox{$\Bcirc$}
  \kern-28pt\rlap{\lower.4pt\hbox{$4$}}\kern 28pt
  \kern-15.5pt\rlap{\raise10pt\hbox{$\node\kern 1pt\scriptstyle t$}}
                 \lower10pt\hbox{$\node\kern 1pt\scriptstyle u$}
 &\BQ(\sqrt{-2})&\cc\mapsto(t,s,u)\\
\hline
\vrule width0pt height 1.2cm depth 1.0cm
G_{13}&\stacktwo{
 \kern 15pt\raise5pt\hbox{$\BBcirc$}\kern-49.5pt
 {\scriptstyle s}\kern 2pt\snode\kern8pt\raise5pt\hbox{$\sbcirc$}
 \kern-11.8pt\rlap{\raise12pt\hbox{$\snode\kern 3pt\scriptstyle t$}}
                 \lower5.5pt\hbox{$\snode\kern 3pt\scriptstyle u$}
 \kern-20pt{\lower 0.5pt\hbox{$\scriptstyle 5$}}
 \kern7pt{\raise2.5pt\hbox{$\scriptstyle 4$}}}
 {\lower10pt\hbox{$\stacktwo{tust=ustu,}{stust=ustus}$}}
 &\BQ(\zeta_8)&
 \stacktwo{\cc\mapsto(s,u,t^s),}{\gal 3\mapsto(s,u^{tu},u^{su})}\\
\hline
\vrule width0pt height 0.8cm depth 0.5cm
G_{14}&\nnode s\overmark8{\bar20pt}\ncnode3t&\BQ(\zeta_3,\sqrt{-2})&
   \stacktwo{\cc\mapsto(s\inv,(t\inv)^s),}{\gal 7\mapsto(s,t^{sttst})}\\
\hline
\vrule width0pt height 1.2cm depth 0.9cm
G_{15}&\stacktwo{
   {\scriptstyle u}\node \kern 0pt
               \rlap{\kern -4pt \lower 10pt\hbox{$\scriptstyle 5$}}
\raise 16.6pt\hbox{$\ucirc$}
               \kern -28.8pt\lower 12.2pt\hbox{$\lcirc$}
\kern -28.8pt \phantom{\dbar 14 pt}
    \kern-2.6pt
    \raise10pt\hbox{$\node$\rlap{\hbox{$\kern 2pt\scriptstyle s$}}}
    \kern-8.8pt
    \lower10.2pt\hbox{$\cnode 3$\rlap{\hbox{$\kern 2pt\scriptstyle t$}}}}
    {\lower-20pt\hbox{$\stacktwo{stu=ust,}{tusts=ustst}$}}
    &\BQ(\zeta_{24})&
    \stackthree{\cc\mapsto(s,t\inv,u^s),}
    {\gal {13}\mapsto(s^{usu},t,u),}
    {\gal {19}\mapsto(s,t,u^{sus})}\\
\hline
\vrule width0pt height 0.5cm depth 0.5cm
G_{16}&\ncnode 5s \bar20pt\ncnode 5t&\BQ(\zeta_5)&
       \gal 2\mapsto(s^2,(t^2)^{s^2t\inv s})\\
\hline
\vrule width0pt height 0.8cm depth 0.5cm
G_{17}&\nnode s\tbar20pt\ncnode 5t&\BQ(\zeta_{20})&
   \stacktwo{\cc\mapsto(s\inv,t\inv),}{\gal 7\mapsto(s^{t\inv st^2},t^2)}\\
\hline
\vrule width0pt height 0.8cm depth 0.5cm
G_{18}&\ncnode 3s\dbar20pt\ncnode 5t&\BQ(\zeta_{15})&
   \stacktwo{\cc\mapsto(s\inv,t\inv),}{\gal 7\mapsto(s^{t\inv s},t^2)}\\
\hline
\vrule width0pt height 0.85cm depth 0.85cm
G_{19}&
 {\scriptstyle s}\node\kern13.5pt\raise2.5pt\hbox{$\Bcirc$}
 \kern-15.5pt\rlap{\raise10pt\hbox{$\cnode 3\kern 1pt\scriptstyle t$}}
 \lower10pt\hbox{$\cnode 5\kern 1pt\scriptstyle u$}
 &\BQ(\zeta_{60})& \stackthree{\cc\mapsto(s^t,t\inv,u\inv),}
   {\gal 7\mapsto(s^{ts},t^{st\inv},u^2),}
   {\gal {41}\mapsto(s,(t\inv)^{(ut)^{su}},u)}\\
\hline
\vrule width0pt height 0.8cm depth 0.5cm
G_{20}&\ncnode3s\overmark5{\bar20pt}\ncnode3t&\BQ(\zeta_3,\sqrt 5)&
   \stacktwo{\cc\mapsto(s\inv,t\inv),}{\gal 7\mapsto(s,s^{(ts\inv)^2})}\\
\hline
\vrule width0pt height 0.85cm depth 0.85cm
G_{21}&\nnode s\overmark{10}{\bar20pt}\ncnode3t&\BQ(\zeta_{12},\sqrt5)&
   \stackthree{\cc\mapsto(s\inv,t\inv),}{\gal 7\mapsto(s^{(t\inv)^{st}},t),}
{\gal {13}\mapsto(s^{t\inv s^{ts}},t)}\\
\hline
\vrule width0pt height 0.8cm depth 0.7cm
G_{22}&
  {\scriptstyle s}\node\kern13.5pt \raise2.5pt\hbox{$\Bcirc$}
  \kern-28pt\rlap{\lower.4pt\hbox{$5$}}\kern 28pt
  \kern-15.5pt\rlap{\raise10pt\hbox{$\node\kern 1pt\scriptstyle t$}}
  \lower10pt\hbox{$\node\kern 1pt\scriptstyle u$}
  &\BQ(\ii,\sqrt 5)&
  \stacktwo{\cc\mapsto(u,t,s),}{\gal 7\mapsto(u^{st},u^{ts},t^{us})}\\
\hline
\vrule width0pt height 0.5cm depth 0.5cm
G_{23}=H_3&\nnode s\overmark5{\bar16pt}\nnode t\bar16pt\nnode u
           &\BQ(\sqrt 5)&\gal 2\mapsto(u^{tsts},u,t)\\
\hline
\vrule width0pt height 1cm depth 0.6cm
G_{24}&
 \nnode s\rlap{\kern 2.5pt{\raise 6pt\hbox{$\triangle$}}}
 \dbar14pt\nnode t
 \kern-26.2pt\raise8.5pt\hbox{$\diagup$}
 \kern -3.1pt
 \raise16.7pt\hbox{$\node$\rlap{\raise 2pt\hbox{$\kern 1pt\scriptstyle u$}}}
 \kern -3pt\raise8.3pt\hbox{$\diagdown$}
 \kern -7.6pt\raise9pt\hbox{$\diagdown$}
 \kern 4pt&\BQ(\sqrt{-7})&\cc\mapsto(u,t,s)\\
\hline
\vrule width0pt height 0.5cm depth 0.5cm
G_{25}&\ncnode3s\bar16pt\ncnode3t\bar16pt\ncnode3u&\BQ(\zeta_3)&
       \cc\mapsto(s\inv,t\inv,u\inv)\\
\hline
\vrule width0pt height 0.5cm depth 0.5cm
G_{26}&\nnode s\dbar16pt\ncnode3t\bar16pt\ncnode3u&\BQ(\zeta_3)&
       \cc\mapsto(s\inv,t\inv,u\inv)\\
\hline
\vrule width0pt height 1cm depth 0.6cm
G_{27}&
 \nnode s\rlap{\kern 2.5pt{\raise 6pt\hbox{$\triangle$}}}
 \dbar14pt\nnode t
 \rlap{\kern-6pt\raise11pt\hbox{$\scriptstyle 5$}}
 \kern-27pt\raise9pt\hbox{$\diagup$}
 \kern -2pt
 \raise16.5pt\hbox{$\node$\rlap{\raise 2pt\hbox{$\kern 1pt\scriptstyle u$}}}
 \kern -1pt\raise9pt\hbox{$\diagdown$}
 \kern 4pt&\stacktwo{\BQ(\zeta_3,\sqrt 5)}{K'=\BQ(\zeta_{15})}&
   \stacktwo{\cc\mapsto(s^{tut},t,u),}{\gal {7}\mapsto(u^{ts},t^u,t)}\\
\hline
\vrule width0pt height 1cm depth 0.5cm
G_{29}&
 \nnode s\bar10pt
 \nnode t\rlap{\kern 1pt\raise9pt\hbox{$\underleftarrow{}$}}
        \rlap{\kern 6.5pt\vrule width1pt height16pt depth-11pt
\kern 1pt \vrule width1pt height16pt depth-11pt}
\dbar14pt\nnode u
 \kern-27pt\raise9.5pt\hbox{$\diagup$}
 \kern -3pt\raise18.5pt
\hbox{$\node$\rlap{\raise 2pt\hbox{$\kern 1pt\scriptstyle v$}}}
 \kern -2pt\raise9.5pt\hbox{$\diagdown$}
 &\BQ(\ii)&\cc\mapsto(s,t,u,v^u)\\
\hline
\vrule width0pt height 0.5cm depth 0.5cm
G_{30}=H_4&\nnode s\overmark 5{\bar10pt}\nnode t\bar10pt\nnode u\bar10pt%
      \nnode v&\BQ(\sqrt5)&\gal 2\mapsto(u^{tsts},u,t,v^{utstuststuvuts})\\
\hline
\vrule width0pt height 1.2cm depth 0.4cm
G_{31}&
   \nnode v\kern -2pt\raise9pt\hbox{$\diagup$}
      \kern -3.5pt\raise18.5pt\hbox{$\nnode s$}
      \kern-13.5pt\bar 19pt\nnode t
   \kern-4pt\raise18.3pt\hbox{$\sBcirc$}
   \kern-17.6pt\bar19pt\nnode w
   \kern-22pt\raise18.3pt\hbox{$\nnode u$}
   \kern-1pt\raise9.8pt\hbox{$\diagdown$}
   &\BQ(\ii)&\cc\mapsto(u,t,s,w,v)\\
\hline
\vrule width0pt height 0.5cm depth 0.5cm
G_{32}&\ncnode3s\bar10pt\ncnode3t\bar10pt\ncnode3u\bar10pt\ncnode3v
       &\BQ(\zeta_3)&\cc\mapsto(s\inv,t\inv,u\inv,v\inv)\\
\hline
\vrule width0pt height 1.1cm depth 0.4cm
G_{33}&\nnode s\bar10pt\trianglerel tvu\bar10pt\nnode w&\BQ(\zeta_3)&
       \cc\mapsto(w,v,u,t,s)\\
\hline
\vrule width0pt height 1.1cm depth 0.4cm
G_{34}&\nnode s\bar10pt\trianglerel tvu\bar10pt\nnode w\bar10pt\nnode x
       &\BQ(\zeta_3)&\cc\mapsto(w,v,u,t,s,x^{wvtuvwstuvts})\\
\end{longtable}

We  have some  observations to  make on  table \ref{titilde}. First, the morphism
$\itilde$  is  not  unique;  often,  there  are  even many different orbits of
possible  $\itilde$ up to inner automorphisms of $G$. Each time we have chosen
one which preserves a parabolic subgroup, which will be used to construct nice
models of the reflection representation of $G$ in section \ref{models}.

This  is not always the choice with  the simplest formulas. In particular, for
Shephard groups (those which have the same diagram, apart from the order of the
reflections,  as a  Coxeter group,  thus the  same braid  group as  a Coxeter
group),  it turns out that it is always  possible to find a morphism which maps
the  complex conjugation $\cc$ to the automorphism which sends each element of
$S$  to  its  inverse.  When $\cc$ does not generate $\Gal(K/\BQ)$,
one can send $\gal{7}$ to $(s^{ts\inv},t\inv)$ for $G_{10}$ and
$\gal{2}$ to $(s^2,(s^2)^{t^{-2} s^2})$ for $G_{16}$ to get such
a morphism (in the case of $G_{14}$, $\gal{7}$ can keep the same image as in 
table \ref{titilde}).
A  general  explanation  for  this phenomenon would be
desirable.

\section{Theorems \ref{A}, \ref{xx} and \ref{C} for the infinite series}
The field of definition of $G(de,e,r)$ is $K=\BQ(\zeta_{de})$, except when $d=
1$  and $r = 2$. In this last  case, $G(e,e,2)$ is the dihedral group of order
$2e$,  for which $K=\BQ(\cos(\frac{2 \pi}{e}))  = \BQ( \zeta_e + \zeta_e\inv)$
and we set $K'=\BQ(\zeta_e)$.

In any case the natural representation of $G(de,e,r)$ via monomial matrices is
over  the vector space $V=K'^r$, and writing $V=V_0\otimes_\BQ K'$ where $V_0$
is  defined by  the standard  basis of  $V$, the  set $G(de,e,r)$  is globally
invariant  under  the  induced  action  of  $\Gal(K'/\BQ)$  on $\GL_r(K')$. It
follows  that this action induces a morphism $\eta : \Gal(\BQ(\zeta_{de})/\BQ)
\to  \Aut(G)$. The purpose of this section  is to show that a slight variation
on  $\eta$ satisfies  the properties  of theorem  \ref{xx}, which  implies the
weaker  theorem  \ref{A}.  Actually  we  will  get  a  more precise version of
\ref{xx},  stating the $\itilde_\chi$ equivariance of some models, which gives
also theorem \ref{C}.

\subsection{Dihedral groups}

In  this  section  we  write  $s'=s'_1$,  $s=s_1$,  thus  the  dihedral  group
$G(e,e,2)$  is generated  by $S=\{s',s\}$  where $$  s' =  \begin{pmatrix} 0 &
\zeta_e\\ \zeta_e\inv & 0 \end{pmatrix}, \quad s = \begin{pmatrix} 0 & 1 \\ 1
&  0 \end{pmatrix} $$ with relations $s^2 =  s^{\prime 2} = 1$, $(s's)^e = 1$.
For  $i$ prime to  $e$ we have  $\eta(\gal i)(s) =  s$ and $\eta(\gal i)(s') =
(s's)^i  s$. In  particular $\eta(\gal  i)(s's) =  (s's)^i$. We may also check
that for any $x\in G$ we have $\eta(\gal{-i})(x) = s \eta(\gal i)(x)s$.

Recall that $K=\BQ(\zeta_e+\zeta_e\inv)$ and let $K'=\BQ(\zeta_e)$. It follows
that  the quotient $\Gamma'=\Gal(K'/\BQ)\to\Gamma=\Gal(K/\BQ)$  is obtained by
identifying $\gal i$ and $\gal{-i}$, thus the above computations show that the
composed  morphism $\Gamma'\xrightarrow\eta\Aut(G)\to\Out(G)$ factors  through $\Gamma$,
leading  to the following commutative diagram  where we have written $\itilde$
for $\eta$.
$$
\xymatrix{
\Gamma' \ar@{^{(}->}_{\itilde}[r] \ar[dr] \ar@{->>}[d] & \Aut(G) \ar@{->>}[d] \\
\Gamma    \ar@{.>}_{\ibar}[r] & \Out(G)
}
$$
The maps $\itilde$ and $\ibar$ above satisfy the hypotheses of
theorems \ref{A} and \ref{xx}.

We    note   that   the   above    reflection   representation   is   actually
$\itilde$-equivariant,  thus  we  have  property  \ref{C}  for  the reflection
representation.

It  is clear that \ref{A}, \ref{xx} and \ref{C} hold for the linear characters
of $G$,
because  they take values in  $\BQ$. The remaining irreducible representations
of $G$ have a model of the form
$$
s' \mapsto \begin{pmatrix} 0 & \zeta_e^r\\ \zeta_e^{-r} & 0 \end{pmatrix},
\quad
s \mapsto \begin{pmatrix} 0 & 1 \\ 1 & 0 \end{pmatrix},
\quad
s's \mapsto \begin{pmatrix} \zeta_e^r & 0  \\  0 & \zeta_e^{-r} \end{pmatrix}
$$
for some $r$, and are determined up to isomorphism by the value of
$\zeta_e^{r}+\zeta_e^{-r}$, so
\ref{A}, \ref{xx} and \ref{C} also hold for them.

\subsection*{Construction of $\iota$.}
We will now determine when it is possible to take $K'=K$.
It is classical that the short exact sequence
$$
1 \to \Gal(\BQ(\zeta_e)/\BQ(\zeta_e+\zeta_e\inv)) \to \Gal(\BQ(\zeta_e) / \BQ) \to
\Gal(\BQ(\zeta_e+\zeta_e\inv) / \BQ) \to 1
$$
splits precisely when $\BQ(\zeta_e)$ contains a quadratic imaginary
extension of $\BQ$. This is the case if and only if $4 | e$ or
$e$ has a prime factor congruent to $3$ modulo $4$.
In that case one can define injective morphisms $\Gamma \to
\Aut(G)$ by composition. 
$$
\xymatrix{
\Gamma' \ar@{^{(}->}_{\itilde}[r] \ar@{->>}[d] & \Aut(G) \ar@{->>}[d] \\
\Gamma  \ar@(ul,dl)@{.>}[u] \ar@{.>}_{\iota}[ur] \ar@{.>}_{\ibar}[r] & \Out(G)
}
$$
We cannot expect to lift this morphism to $\Aut(G)$ in general, as already
shows the example $e=5$. In that case, it is easily checked
that $\Gamma = \Gal(\BQ(\sqrt{5})/\BQ)$ has order two and
that all involutive automorphisms of $G$ are inner.

If the exact sequence splits, we let $H < \Gamma'$ be a section
of $\Gamma$. Since $H$ has index two in $\Gamma'$ we know that
$H$ is normal in $\Gamma'$ and, by elementary Galois theory
we deduce from $K = K' \cap \BR$ that $H = \Gal(K'/\BQ(\ii \alpha))$
for some quadratic number $\alpha \in \BR$. Let $\rho : G \to \GL_2(K')$
denote an arbitrary 2-dimensional irreducible representation
of $G$, following the models described above. An intertwiner
between $\rho$ and $\cc \circ \rho$ is given by $\rho(s) \in \GL_2(\BQ)$.
Let
$$
M = \frac{1}{4\ii \alpha} \begin{pmatrix}1+\ii \alpha&-1+\ii \alpha\\-1+\ii \alpha&1+\ii \alpha\end{pmatrix}
\in \GL_2(\BQ(\ii \alpha)).
$$
It is easily checked that $M\inv \cc(M)$ is equal to $\rho(s)$ up to a scalar
multiple (a general way to find such an $M$ will be described in section
\ref{galoissect}). Let $\rho' : G \to \GL_2(K')$ be defined by
$\rho'(g) = M \rho(g) M\inv$. We have $\cc(\rho'(g)) = 
M \rho(s) \cc(\rho(g)) \rho(s) M \inv = M \rho(g) M \inv = \rho'(g)$
hence $\rho'(g) \in \GL_2(K)$. Moreover, for $\gamma \in H$ we
have $\gamma(M) = M$ hence $\gamma \circ \rho'(g) = M \gamma \circ \rho(g) M \inv
= M \circ \rho \circ \itilde(\gamma) (g) M \inv = \rho' \circ \itilde(\gamma) (g)$.
It follows that, in these cases, we have morphisms $\iota : \Gamma \into \Aut(G)$
and models over $K$ for all irreducible
representations $\rho$ of $G$, such that $\gamma \circ \rho = \rho \circ \iota(\gamma)$
for all $\gamma \in \Gamma$. Thus, in those cases we can take $K' = K$ in
theorems \ref{xx} and \ref{C}.

\subsection{The case $G(d,1,r)$}

The   group  $G=G(d,1,r)$  is  generated  by  the  set  of  reflections  $S=\{
t,s_1,\ldots,s_{r-1}\}$  where $t=\Diag(\zeta_d,1,\ldots,1)$
and $s_k=M_{(k,k+1)}$ with the notations of the proof of \ref{refl}.  We have
$\eta(\gal{\alpha})(t) = t^{\alpha}$ and $\eta(\gal{\alpha})(s_k) = s_k$.

We  provide an easy  proof of \ref{A}  and \ref{xx} for  $G(d,1,r)$ starting
from  the  case  $r  =  2$.  Then  $G$  has  order $2d^2$, and its irreducible
representations  have dimension 1 and  2. The case of  linear characters is
immediate to check.   The   2-dimensional   representations   are   given,  for
$\zeta_1,\zeta_2 \in \mu_d, \zeta_1\ne\zeta_2$, by the model
$$
\rho_{\zeta_1,\zeta_2}(t)=\begin{pmatrix}\zeta_1&0\\0&\zeta_2\end{pmatrix}
\qquad
\rho_{\zeta_1,\zeta_2}(s_1)=\begin{pmatrix}0 & 1 \\ 1 & 0 \end{pmatrix}
$$
It is readily checked that
$$
\rho_{\zeta_1,\zeta_2} \circ \eta(\gal{\alpha}) = \gal{\alpha} \circ
\rho_{\zeta_1,\zeta_2} = \rho_{\zeta_1^{\alpha},\zeta_2^{\alpha}}
$$
and  this proves theorems \ref{A}, \ref{xx} for the groups $G(d,1,2)$ by taking
$K'=K$  and $\itilde=\eta$. The first equality in the above formula
also shows
that the model is $\eta$-equivariant, thus giving \ref{C}.

Theorems \ref{A} and \ref{xx} are then easily deduced from this for the groups
$G(d,1,r)$, using that, for $r  \geq 3$, each irreducible
representation   of  $G(d,1,r)$  is  determined   up  to  isomorphism  by  its
restriction    to    the    parabolic    subgroup   $G(d,1,r-1)$   (see   e.g.
\cite{ariki-koike} corollary 3.12).

We  will  give  a  proof  of  \ref{C}  for  $G(d,1,r)$  in  the same spirit in
\ref{rational  for primitive}. But in the next subsections we will more
constructively show that the models of \ref{RepGd1r} and \ref{impr}
are $\itilde_\chi$-equivariant for a suitable $\itilde_\chi$.

\subsection{The case of $G(d,1,r)$}

From the formulas $\eta(\gal{\alpha})(t) = t^{\alpha}$ and
$\eta(\gal{\alpha})(s_k)  = s_k$ and the formulae of subsection \ref{RepGd1r},
it  is clear that our models over $\BQ(\zeta_{d})$ satisfy theorem \ref{xx} by
taking  $\itilde=\eta$; actually, the formulae show that the model $\rho_\bla$
is $\itilde$-equivariant, which gives theorem \ref{C}.

More precisely, if we define the operator
$\Sigma_\alpha:\CT(\bla)\to\CT(\gal\alpha(\bla))$ by
$\Sigma_\alpha(T)=(T_0,T_\alpha,T_{2\alpha},\ldots, T_{(d-1)\alpha})$, we
have $\Sigma_\alpha\circ\gal\alpha\circ\rho_\bla=
\rho_{\gal\alpha(\bla)}\circ\Sigma_\alpha$ where 
$\rho_{\gal\alpha(\bla)}$ is as in \ref{gal on repGd1r}.

\subsection{The general case $G(de,e,r)$}
Using that the matrix $S$ of subsection \ref{impr} is rational in
our chosen basis of $V(\bla)$, we get for $\omega\in\mu_{\#\centl}$:
$$\gal\alpha(\chi_{\bla,\omega}(g))=\frac 1{\#\centl}\sum_{i=0}^{\#\centl
-1}\gal\alpha(\omega)^{-i}\Trace(\gal\alpha(\rho_\bla(g))S^i).$$
From this we get, if we write $S_\bla$ for $S$ to keep track where it acts:
$$\gal\alpha(\chi_{\bla,\omega}(g))=\frac 1{\#\centl}\sum_{i=0}^{\#\centl
-1}\omega^{-i\alpha}\Trace(\Sigma_\alpha\inv\rho_{\gal\alpha(\bla)}(g)
\Sigma_\alpha S_\bla^i).$$
Now, it is easy to check that $\Sigma_\alpha S_\bla\Sigma_\alpha\inv=
S_{\gal\alpha(\bla)}^\alpha$, so
$$\gal\alpha(\chi_{\bla,\omega}(g))=\frac 1{\#\centl}\sum_{i=0}^{\#\centl
-1}\omega^{-i\alpha}\Trace(\rho_{\gal\alpha(\bla)}(g)
S_{\gal\alpha(\bla)}^{i\alpha})=\chi_{\gal\alpha(\bla),\omega}(g).$$
While we have by a similar computation
$$\begin{aligned}
\chi_{\bla,\omega}(\itilde(\gal\alpha)(g))&=
\frac 1{\#\centl}\sum_{i=0}^{\#\centl-1}\omega^{-i}\Trace(
\gal\alpha(\rho_\bla(g))S_\bla^i)\\
&=\frac 1{\#\centl}\sum_{i=0}^{\#\centl
-1}\omega^{-i}\Trace(\rho_{\gal\alpha(\bla)}(g)
S_{\gal\alpha(\bla)}^{i\alpha})\\
&=\chi_{\gal\alpha(\bla),\gal\alpha\inv(\omega)}(g).
\end{aligned}$$

We  need slightly more precise formulas that we will get using
our model which simplifies Ariki's one. We claim that if $\theta=\zeta_e^b$
then
$$\gal\alpha\circ\rho_{\bla,\theta^i}=
\rho_{\bla,\theta^i}\circ\Ad t^{i(1-\alpha)}\circ\eta(\gal\alpha)=
\rho_{\bla,\theta^i}\circ\Ad t^i\circ\eta(\gal\alpha)\circ\Ad t^{-i}
$$
indeed this  formula is obvious for the generators $t$ and $s_i$ for $i\ne 1$. 
It is  thus sufficient to check  it for $s_1^{t^k}$; one  uses the
formulae \ref{rhobla1}
and  the fact that $\rho_{\bla,1}(s_1)$ is a  rational matrix. 

Let us define  $\itilde_{\chi_{\bla,\theta^i}}:\Gamma\to\Aut(G)$          by
$\itilde_{\chi_{\bla,\theta^i}}(\gal\alpha)=\Ad
t^i\circ\eta(\gal\alpha)\circ\Ad    t^{-i}$.   The   formula   above   becomes
$\gal\alpha\circ\rho_{\bla,\theta^i}=\rho_{\bla,\theta^i}\circ
\itilde_{\chi_{\bla,\theta^i}}(\gal\alpha)$, \ie. that
$\rho_{\bla,\theta^i}$ is $\itilde_{\chi_{\bla,\theta^i}}$-equivariant,
which proves theorems \ref{xx} and \ref{C}.
We  note that the various morphisms
$\itilde_{\chi_{\bla,\theta^i}}$  differ  by conjugacy by an
element  of $N$. Indeed they are  conjugate by a power of  $\Ad t$ which is a
generator of $N$. When $\gcd(e,r)=1$ then $\Nbar$ is trivial, and
when $\gcd(e,r)=2$ then $\Nbar$ is of order 2 and it results that 
$\itilde_{\chi_{\bla,\theta^i}}$ does not depend on $i$: indeed $\alpha$ being
prime to $de$ is odd thus $1-\alpha$ is even and $\Ad t^{i(1-\alpha)}$ is
always inner, which takes care of the remarks in \ref{A}.

\section{\label{galoissect}Galois descent}

The  existence of a model over  $K$ globally invariant by $\Gamma=\Gal(K/\BQ)$
of  a faithful  representation of  $G$, namely  its reflection representation,
enabled   us   in   the   $G(de,e,r)$   case   to   construct   the   morphism
$\itilde : \Gamma\to\Aut(G)$. 

In  this section,  we show  how by  Galois descent,  one can conversely try to
construct a globally  invariant  model for any given irreducible representation
--- more precisely we will
try to get a $\itilde$-equivariant model. The possible obstruction to do so is
an  element  of  the  Brauer  group  of  $K'$  that turns out to be always
trivial except for the group $G_{22}$. 

\subsection{Non-abelian Galois cohomology}

Recall from \cite{SERRE} that, if $H$ is a group acting on the left on a group
$Q$, then a map $h \mapsto A_h$ from $H$ to $Q$ is a \emph{cocycle} if $A_{st}
= A_s s(A_t)$ for all $s,t \in H$. A particular cocycle, called the \emph{zero
cocycle},  is given by  sending all elements  of $H$ to  the trivial element of
$Q$. The set of cocycles is denoted $Z^1(H,Q)$. Two cocycles $h \mapsto A_{h}$
and  $h \mapsto B_{h}$ are  said to be cohomologous  if there exists $a \in Q$
such  that $B_h =  a\inv A_h h(a)$.  This defines an  equivalence relation on
$Z^1(H,Q)$  whose set of equivalence classes is denoted $H^1(H,Q)$. The set of
cocycles  cohomologous to the  zero cocycle is  called the set of coboundaries
and  is denoted $B^1(H,Q)$.  If $Q$ is  commutative these definitions coincide
with those of ordinary group cohomology.

\subsection*{Application to \ref{C}}

Let  $\rho:G\to\GL(E)$  where  $E$  is  a  $K'$-vector space be an irreducible
representation.  Assume  we  choose  a  $\BQ$-form $E=E_0\otimes_\BQ K'$. This
defines  an  action  of  $\Gamma'=\Gal(K'/\BQ)$  on  $E$  and  by  \ref{A} for
$\gamma\in\Gamma'$  we  have  $\gamma\circ\rho\simeq\rho\circ\itilde(\gamma)$.
This  means that that there exists $A_\gamma\in\GL(E)$ such that for any $g\in
G$   we  have  $A_\gamma\gamma(\rho(g))A_\gamma\inv=\rho(\itilde(\gamma)(g))$.
Since  $\rho$ is absolutely irreducible, by  Schur's lemma, $A_\gamma$ gives a
well-defined   element   of   $\PGL(E)$.   It   is  immediate  to  check  that
$\gamma\mapsto  A_\gamma$  is  in  a  fact  a  cocycle, that is, an element of
$Z^1(\Gamma',\PGL(E))$.

Assume  now that $\rho$ has a globally  invariant model, of the form $g\mapsto
a\rho(g)   a\inv$;   this   means   that   there   exists   a   map   $\iota':
\Gamma'\to\Aut(G)$  such  that $\gamma(a\rho(g)a\inv)=a\rho(\iota'(\gamma)(g))
a\inv$, or in other terms that the model is $\iota'$-equivariant
for some $\iota':
\Gamma'\to\Aut(G)$.  If  $\iota'=\itilde$   we  get   the  equality   in  $\PGL(E)$  that
$a\inv\gamma(a)=A_\gamma$,  that  is  that $\{ \gamma \mapsto A_\gamma \}\in B^1(\Gamma',\PGL(E))$.
Thus  the obstruction to the existence of $\itilde$-equivariant model is an
element of $H^1(\Gamma',\PGL(E))$.

In  practice to apply this procedure  we encounter another problem which turns
out  to involve also Galois  cohomology: although we know that 
any $\rho$ has  a model over  $K'$ (even $K$), sometimes
(\eg.  when using the \CHEVIE\ data) we are not given such a model but a model
$E=E_0\otimes_{K'}L$  over some Galois  extension $L$ of  $K'$. Again the fact
that  for $\gamma\in\Pi=\Gal(L/K')$ we  have $\gamma\circ\rho\simeq\rho$ means
that  we can find an intertwiner $A_\gamma$ such that $A_\gamma\gamma(\rho(g))
A_\gamma\inv=\rho(g)$,  which is  an element  of $Z^1(\Pi,\PGL(E))$; and again
the  existence of  a model  over $K'$  is equivalent  to this  cocycle to be a
coboundary.

We will now give an algorithm to check the vanishing of a cocycle of a Galois
group into the projective linear group.

\subsection{Brauer groups}
 
In  this subsection and the next one, $K_0 \subseteq K$ will denote an arbitrary
Galois  extension  of  number  fields.  Let  $\Gamma  =  \Gal(K/K_0)$  and let
$\Br(K/K_0)  =  H^2(\Gamma,K^\times)$  be  the  Brauer  group.  Let  $E  = E_0
\otimes_{K_0} K$ be a $K_0$-form of a finite dimensional $K$-vector space $E$.
The  short exact sequence  $1 \to K^\times  \to \GL(E) \to  \PGL(E) \to 1$ gives
rise  to  a  Galois  cohomology  long  exact  sequence  (of  pointed sets; see
\cite[prop. 2 p. 133]{SERRE}):
\begin{multline*}
1\to H^0(\Gamma,K^\times)\to H^0(\Gamma,\GL(E)) \to H^0(\Gamma,\PGL(E))
\to H^1(\Gamma,K^\times)\to\\ 
\hfill H^1(\Gamma,\GL(E))\to H^1(\Gamma,\PGL(E))\to H^2(\Gamma,K^\times)\\
\end{multline*}
and    in    particular    a    coboundary   operator   $H^1(\Gamma,\PGL(E))\to
H^2(\Gamma,K^\times)$.   This  map  sends  $\{\gamma  \mapsto  A_{\gamma}\}\in
H^1(\Gamma,\PGL(E))$  to the class of $(\gamma,\tau) \mapsto \tilde{A}_{\gamma}
\gamma(\tilde{A}_{\tau})        \tilde{A}_{\gamma        \tau}\inv$        in
$H^2(\Gamma,K^\times)$,     where,    for    all    $\gamma    \in    \Gamma$,
$\tilde{A}_{\gamma}$  is  a  representative  of  $A_{\gamma}$ in $\GL(E)$. This
definition   does   not   depend   on   the   choice  of  the  representatives
$\tilde{A}_{\gamma}$  and this  map is  known to  be injective  (this uses the
version  of  Hilbert's  theorem  90  which  says  that $H^1(\Gamma,\GL(E))$ is
trivial; see \cite{SERRE} ch. X proposition 8 and proposition 9).

In  the important case where $K$ is a cyclic extension of degree $n$ of $K_0$,
assume  $\Gamma=\langle\gamma\rangle$ and let $N : K  \to K_0$ be the norm map
$x\mapsto   x\gamma(x)\ldots\gamma^{n-1}(x)$.  Then  an  explicit  isomorphism
from $H^2(\Gamma,K^\times)$ to the $0$-th Tate cohomology group $\hat{H}^0   (\Gamma,K^\times)  =  K_0^\times/N
K^\times$ is given by the Nakayama map
$$  c  \mapsto \prod_{k=0}^{n-1} c(\gamma^k,  \gamma)$$ 
hence    $\{\gamma\mapsto   A_{\gamma}\}$   is   sent    to   the   class   of
$\tilde{A}_{\gamma}      \gamma(\tilde{A}_{\gamma})     \dots     \gamma^{n-1}
(\tilde{A}_{\gamma})$   (a  scalar   matrix  identified   to  an   element  of
$K_0^\times$) in $K_0^\times/N K^\times$.

It  follows that, when $K$ is a cyclic extension of $K_0$, the vanishing of $c
\in  H^2(\Gamma,K^\times) = \Br(K/K_0)$ boils down to the verification
that  some element of $K_0^\times$ is a norm of an element of $K$. When $K$ is
an  \emph{abelian} extension of  $K_0$, we can  always reduce to  this case by
induction  on  $[K:  K_0]$,  by  the  inflation-restriction exact sequence for
Brauer  groups.  Indeed,  choosing  a  cyclic  subgroup  $\Gamma'$ of $\Gamma$
defines a Galois sub-extension $K_1 = K^{\Gamma'}$ of $K$ such that
$\Gal(K/K_1)$ is cyclic. We then have
the exact sequence (see \cite{SERRE} ch. X \S 4 proposition 6)
$$
0 \to \Br(K_1/K_0) \to \Br(K/K_0) \to \Br(K/K_1)
$$
and  we may  follow the  following procedure  : first  check that the image in
$\Br(K/K_1)$ is zero, and then start again with the induced element in
$\Br(K_1/K_0)$, until $K_1/K_0$ itself becomes a cyclic extension.

We finally mention the following general result.
\begin{proposition}  \label{odd  dimension}  If  $\# \Gamma$  is prime
to $\dim(E)$, then $H^1(\Gamma,\PGL(E))=0$
\end{proposition}
\begin{proof}
Let $N = \dim(E)$. Let  $c\in  H^1(\Gamma,\PGL(E))$  and  let  $A$  be  the central simple
algebra associated  to it. Let $e$ be its exponent, that is
the order of its class $[A]$ in the Brauer group of $K_0$.
Since $A \otimes K \simeq M_N(K)$ we have $[A] [K] = 0$ in the
Brauer group hence $e| [K:K_0]$ and $e | \# \Gamma$ (see \cite{SERRE}
X \S 4 ex. 2). On the other hand, $A \simeq M_r(D)$ with $D$ a skew field
whose center is $K_0$ and $\dim_{K_0} A = N^2$. It follows that
$N^2 = [A : K_0] = r^2 [D: K_0]$. Now $[D:K_0] = m^2$ for some integer
$m$ (called the index of $D$). It is a classical fact that the exponent
divides the index (\cite{SERRE} X \S 5 ex. 3a), thus $e | N$,
and $e = 1$ since $N$ and $\# \Gamma$ are coprime. The conclusion
follows from the fact that  $H^1(\Gamma,\PGL(E))$ embeds in the Brauer group
of $K_0$.
\end{proof}

\subsection*{Explicit computations}

In  practice  (\eg.  to  get  an  explicit  model  over  a  smaller field of a
representation),  it  is  not  enough  to  solve  the problem whether a cocycle
$\{\gamma\mapsto  A_\gamma\}\in Z^1(\Gamma,\PGL(E))$ is  a coboundary. We want
an  explicit element $M \in \PGL(E)$ such that $A_{\gamma} = M\inv \gamma(M)$
for  all  $\gamma\in  \Gamma$.  If  we  can  get  a preimage $(B_{\gamma})$ of
$(A_{\gamma})$  in  $Z^1(\Gamma,\GL(E))$,  then  by  Hilbert's Theorem 90 this
cocycle  is a coboundary, meaning that there exist $\tilde{M} \in \GL(E)$ such
that  $B_{\gamma}  =  \tilde{M}\inv  \gamma(\tilde{M})$  for  all $\gamma \in
\Gamma$.  A rather constructive  proof of  this theorem  (see e.g. \cite{SERRE} ch. X
proposition 3) goes as follows. For all $C \in \End(E)$, the expression
$$
X = \sum_{\gamma \in \Gamma} B_{\gamma} \gamma(C)
$$
satisfies  $\gamma(X)  =  B_{\gamma}\inv  X$.  Because $K$ has characteristic
zero,  general arguments imply that ``many'' $C  \in K$ exist such that $X$ is
invertible.  If this  is the  case, then  $B_{\gamma} = X \gamma(X)\inv$ thus
$\tilde{M} = X\inv$ is the desired solution. For instance, if $\Gal(K / K_0)
=\{1,\gamma\}$
has  order 2 with generator $\gamma_0$, then  the condition on  $C$ is that  $- C \gamma(C)$ is not an
eigenvalue for $B_\Id\inv B_{\gamma_0}$.

In  the  case  where  $K$  is  a  cyclic  extension of $K_0$, the study of the
previous  section shows  that we  can always  lift a  coboundary to  a cocycle
$\{\gamma\mapsto      B_{\gamma}\}$      for      $\GL(E)$.     Indeed,     if
$\Gamma=\langle\gamma\rangle$    we   have   seen   that   $\tilde{A}_{\gamma}
\gamma(\tilde{A}_{\gamma})     \dots     \gamma^{n-1}(\tilde{A}_{\gamma})    =
N(\lambda)\Id$   for  some  $\lambda  \in  K^\times$  so  that  $B_{\gamma}  =
\lambda\inv  \tilde{A}_{\gamma}$ belongs to  $Z^1(\Gamma,\GL(E))$ (it is easy
to  check  that the  equation $\Id=\tilde
A_\gamma\ldots  \gamma^{n-1}(\tilde A_\gamma)$  is necessary  and sufficient
for $\tilde A_{\gamma^i}=\tilde
A_\gamma\ldots  \gamma^{i-1}(\tilde A_\gamma)$ to define
a cocycle of $\langle\gamma\rangle$).

\subsection{An algorithm to check \ref{C}}
Let  $\rho : G  \to \GL(E)$ be  an irreducible representation  of $G$ over the
$K'$-vector  space $E$  for which  we would  like to  get a $\itilde$-equivariant
model.  When $K'/\BQ$ is cyclic, we found an explicit necessary and sufficient
condition  for a  class in  $H^1(\Gamma',\PGL(E))$ to  vanish and thus to give an
equivariant model.

In  general we know that $K'/\BQ$ is an abelian extension. We will show how to
reduce  by induction to the case of a cyclic extension, by making explicit the
inflation-restriction exact sequence in this case.

Choose  a tower of  (Galois) extensions $\BQ  = K_0 \subset  K_1 \subset \dots
\subset  K_m = K'$  such that $K_r/K_{r-1}$  is a cyclic  extension for all $1
\leq  r \leq m$. We fix a $\BQ$-form $E_0$ of $E$ and a basis of $E_0$ so that
we identify $E_0$ with $\BQ^n$ and $E$ with $K^{\prime n}$.

Assume that, for some $r \in [1,m]$ we have a model $\rho_r : G \to \GL(E)$ of
$\rho$  such  that  $\gamma  \circ  \rho  =  \rho \circ \iota(\gamma)$ for all
$\gamma \in \Gal(K'/K_r)$. We will show how to decrease $r$.

Let  $\{\tau\mapsto M_{\tau}\} \in Z^1(\Gal(K'/K_{r-1}),\PGL(E))$ intertwining
$\tau   \circ   \rho_r$   and   $\rho_r   \circ  \iota(\tau)$  for  $\tau  \in
\Gal(K'/K_{r-1})$.  For  $\gamma$  in  the  subgroup  $\Gal(K'/K_r)$  we  have
$M_{\gamma}=1$  by the induction  hypothesis; we thus  have $M_{\tau \sigma} =
M_{\tau}$  for all $\sigma \in \Gal(K'/K_r)$ by the cocycle condition $M_{\tau
\gamma}  = M_{\tau}  \tau(M_{\gamma})$. Hence  $M_{\tau}$ only  depends on the
class  of  $\tau$  in  $\Gal(K_r/K_{r-1})$,  and lies in $\PGL_n(K_r)$ because
$M_{\tau}  = M_{\tau \gamma} = M_{\gamma \tau} = M_{\gamma} \gamma(M_{\tau}) =
\gamma(M_{\tau})$   for  all   $\gamma  \in   \Gal(K'/K_r)$.  It   follow  that
$\{\tau\mapsto       M_{\tau}\}$       defines       an       element       of
$Z^1(\Gal(K_r/K_{r-1}),\PGL_n(K_r))$;  we are  thus reduced  to a cyclic case;
assuming  we can solve it  we find $X \in  \GL_n(K_r)$ such that $M_{\tau} = X
\tau(X)\inv$  for all  $\tau \in  \Gal(K'/K_{r-1})$. Then  $\rho_{r-1}(g) = X
\rho_r(g)  X\inv$ is a model  of $\rho$ such that  $\gamma \circ \rho_{r-1} =
\rho_{r-1} \circ \iota(\gamma)$ for all $\gamma \in \Gal(K'/K_{r-1})$.

We  could successfully  carry out  this algorithm  for all  representations of
exceptional  groups for which we had a model. However, with the exception of the
cyclic  case,  this  procedure  does  not  lead  to necessary conditions. Note
however the following observation :

\begin{proposition} \label{oddinvariant} Assume $G$ to be an exceptional
irreducible complex reflection group. Then any odd-dimensional
irreducible representation has a globally invariant model over $K'$.
\end{proposition}
\begin{proof}
A  case-by-case analysis shows that $[K':\BQ]$ is always a power of $2$. The
result is then an immediate consequence of \ref{odd dimension}.
\end{proof}

\section{\label{models} Invariant models}

We will now prove theorem \ref{C} using the results of the previous sections
and multiplicity one property of tensor products.

For two representations $\rho_1 , \rho_2$ of a group $G$,
we will denote $(\rho_1\mid \rho_2)$ the dimension of the
space of intertwiners of $\rho_1$ and $\rho_2$. This coincides
with the usual scalar product of the corresponding characters.

We will make repeated use of the following lemma.

\begin{lemma} \label{mfree}Let $G$ be a finite group and $K$ be a field
of characteristic 0.
Let $\rho'$ be a finite-dimensional
representation of $G$ defined over $K$
and let $\rho \in \Irr(G)$
be such that $(\rho'\mid \rho) = 1$.
Let $\chi$ be the character of $\rho$.
Assume that $\chi$ takes values in $K$. Then $\rho$ admits a model over $K$
and :
\begin{itemize}
\item if $\gamma \in \Gal(K/\BQ)$ and $a \in \Aut(G)$ satisfy
$\gamma \circ \rho' = \rho' \circ a$ and
$\gamma \circ \chi = \chi \circ a$, then $\rho$ admits a model
over $K$ such that $\gamma \circ \rho = \rho \circ a$.
\item if there exists $j : \Gal(K/\BQ) \to \Aut(G)$ such that
$\gamma \circ \rho' = \rho' \circ j(\gamma)$ and
$\gamma \circ \chi = \chi \circ j(\gamma)$, then $\rho$
admits a $j$-equivariant model over $K$.
\end{itemize}
\end{lemma}
\begin{proof} Since $\rho$ is an isotypic component of
a $\rho'$ and has its character over $K$ it is defined
over $K$ ; indeed, the $G$-equivariant projector
$$
p = \frac{\chi(1)}{\# G} \sum_{g \in G} \chi(g\inv) \rho'(g)
$$
on  the $\rho$-isotypic component of $\rho'$ belongs to $\End_K(V)$. Moreover,
if  $\gamma \circ \chi = \chi \circ a$,  we have $$ \gamma(p) = \frac
{\chi(1)}{\# G}
\sum_{g \in G} \chi(a(g)\inv) \rho'(a(g)) = p $$ thus a $\BQ$-form $V' = V'_0
\otimes_\BQ K$ induces a $\BQ$-form on $V = \Image p$. It follows that $\gamma
\circ \rho = \rho \circ a$, which proves the last two points.
\end{proof}

\subsection{Reflection representations}

We now prove theorem \ref{C} for the reflection representation of exceptional
groups, that is
the following proposition.

\begin{proposition} \label{reflprim}
Let $G$ be an exceptional irreducible complex reflection different
from $G_{22}$ given by
$\rho : G \into \GL(V)$ where $V$ is a $K$-vector space and $K$
is the field of definition of $G$.
Let $K'$ be as in theorem \ref{xx}. Then $\rho$ has a
$\itilde_{\rho}$-equivariant model over $K'$.
\end{proposition}
Recall that $K = K'$ except when $G = G_{27}$.

\subsection*{Models for reflection representations of 2-dimensional groups from the braid group}

A  way  to  obtain  the  2-dimensional  representations of the groups $G_4$ to
$G_{22}$  is by using the  following matrices: 
$$\bs\mapsto
\begin{pmatrix}x_1&\frac{y_1+y_2}{y_1y_2}-\frac{(z_1+z_2)x_2}r\\0&x_2
\end{pmatrix},
\bt\mapsto
\begin{pmatrix} y_1+y_2&1/x_1\\ -y_1y_2x_1&0 \end{pmatrix},
\bu\mapsto
\begin{pmatrix} 0&\frac {-r}{y_1y_2x_1x_2}\\ r&z_1+z_2 \end{pmatrix}$$
where $r=\sqrt{x_1x_2y_1y_2z_1z_2}$.
The braid groups of $G_7$,
$G_{11}$  and $G_{19}$ are isomorphic to the same group $B$, with presentation
$\langle  \bs, \bt, \bu \mid \bs\bt\bu=\bt\bu\bs=\bu\bs\bt \rangle$; the above
matrices  give the  2-dimensional representation  of $B$ where the eigenvalues
of $\bs$ (resp. $\bt$, $\bu$) are $x_1, x_2$ (resp. $y_1,y_2$, $z_1,z_2$).
This  representation factors  through the  Hecke algebra,  the quotient of the
group    algebra   of    $B$   by    the   relations   $(\bs-x_0)(\bs-x_1)=0$,
$(\bt-y_0)(\bt-y_1)(\bt-y_2)=0$   and  $\prod_{i=0}^{i=n-1}(\bu-z_i)=0$  where
$n=3$  (resp. 4,5)  for $G_7$  (resp. $G_{11}$,  $G_{19}$). In  turn the group
algebra of $G_7$ (resp. $G_{11}$, $G_{19}$) is the specialization of the Hecke
algebra  for  $x_i\mapsto  (-1)^i,  y_i\mapsto \zeta_3^i,z_i\mapsto \zeta_3^i$
(resp.  $z_i\mapsto \zeta_4^i$, $z_i\mapsto\zeta_5^i$); the Hecke algebras for
$G_4$  to $G_6$ are  subalgebras of partial  specializations of that for $G_7$
(the same holds for $G_8$ to $G_{15}$ with respect to $G_{11}$ and $G_{16}$ to
$G_{22}$ with respect to $G_{19}$): in each case, these algebras are generated
by  conjugates of  a part  of the  generators (or  of some power of them); the
other generators are specialized to the group algebra.

For $G_4, G_8$ and $G_{16}$ the Hecke algebra is generated by $\bu$ and
$\lexp\bs\bu$.

For $G_5, G_{10}$ and $G_{18}$ it is generated by $\bt$ and $\bu$.

For $G_6, G_9$ and $G_{17}$ it is generated by $\bs$ and $\bu$.

For $G_{14}$ and $G_{21}$ it is generated by $\bs$ and $\bt$.

For $G_{12}$ and $G_{22}$ it is generated by $\bs,\lexp\bt\bs$ and $\bs^\bt$.

For $G_{20}$ it is generated by $\bt$ and $\lexp\bs\bt$.

For $G_{13}$ it is generated by $\bu^2,\bs$ and $\bs^\bt$.

Finally for $G_{15}$ it is generated by $\bs,\bt$ and $\bu^2$.

\subsection*{Proof of proposition \ref{reflprim}} 
The statement is clear if the rank
of $G$ is odd, by proposition \ref{oddinvariant}. Then
we may assume that, either $G = G_{34}$, or that its rank is $2$ or
$4$. Recall that the groups of rank 2 are numbered $G_4,\dots,G_{22}$.

We  started from a model  of the reflection representations
coming  either from the Hecke algebra as  described above, or by root diagrams
as  considered  in  \cite{cohen}. Then we apply the algorithm of Galois
descent; in order to solve norm  equations we used the computer system
\MAGMA. All equations thus obtained for the reflection representation
were solvable except for $G_{22}$.

Note that, even when the original model is globally invariant, the
induced morphism $\eta : \Gamma \to \Aut(G)$ does not in general
coincide with $\itilde_{\rho}$. Furthermore, the induced action of $\eta$
on the isomorphism classes of irreducible representations of $G$
may also differ from the natural Galois action of $\Gamma$,
as shows the second example below.

The  models  given  in  table  \ref{tablemodels} were  simplified  by  the following
elementary  observation. For  all the  images $\tilde{\iota}  (\Gamma)$ chosen
here,  it turns out that  there is a generator  $s$ such that the subgroup $H$
generated  by $s$  is stable  under $\tilde{\iota}(\Gamma)$. Lemma \ref{mfree}
thus  implies that there is  a basis of the  underlying $\BQ$-form of $K^2$ on
which the action of $s$ is diagonal.

\begin{longtable}[c]{|>{$}l<{$}|>{$}l<{$}|}
\caption{Invariant models for $G_4$ to $G_{21}$}\label{tablemodels}\\
\hline G_4&\vrule width0pt height 0.7cm depth 0.5cm 
s\mapsto\begin{pmatrix}1&0\\0&\zeta_3\end{pmatrix},\quad
t\mapsto\frac{1}{\sqrt{-3}}
\begin{pmatrix}-1&\zeta_3\\2&\zeta_3\end{pmatrix}
\\\hline G_5&\vrule width0pt height 0.7cm depth 0.5cm 
s\mapsto\begin{pmatrix}1&0\\0&\zeta_3\end{pmatrix},\quad
t\mapsto\frac1{\sqrt{-3}}
 \begin{pmatrix}\zeta_3&\zeta_3\\2&-1\end{pmatrix}
\\\hline G_6&\vrule width0pt height 0.7cm depth 0.5cm 
s\mapsto\frac1{\sqrt3}\begin{pmatrix}1&1\\2&-1\end{pmatrix},\quad
t\mapsto\begin{pmatrix}1&0\\ 0&\zeta_3\end{pmatrix}
\\\hline G_7&\vrule width0pt height 0.7cm depth 0.5cm 
s\mapsto\frac1{\sqrt3}\begin{pmatrix}1&1\\ 2&-1\end{pmatrix},\quad
t\mapsto\begin{pmatrix}1&0\\0&\zeta_3\end{pmatrix},\quad
u\mapsto\frac1{\sqrt{-3}}
 \begin{pmatrix}\zeta_3&\zeta_3\\2&-1\end{pmatrix}
\\\hline G_8&\vrule width0pt height 0.7cm depth 0.5cm 
s\mapsto\begin{pmatrix}1&0\\0&\ii\end{pmatrix},\quad
t\mapsto\frac{\ii-1}2\begin{pmatrix}-\ii&1\\1&-\ii\end{pmatrix}
\\\hline G_9&\vrule width0pt height 0.7cm depth 0.5cm 
s\mapsto\frac12\begin{pmatrix}\sqrt2&2\\1&-\sqrt2\end{pmatrix},\quad
t\mapsto\begin{pmatrix}1&0\\0&\ii\end{pmatrix}
\\\hline G_{10}&\vrule width0pt height 0.7cm depth 0.5cm 
s\mapsto\frac1{\zeta_3(\ii-1)}\begin{pmatrix}-\ii&2\ii\\ \frac12&1\end{pmatrix},\quad
t\mapsto\begin{pmatrix}1&0\\0&\ii\end{pmatrix}
\\\hline G_{11}&\vrule width0pt height 0.7cm depth 0.5cm 
s\mapsto C=\frac1{\sqrt6}\begin{pmatrix}-2&1\\2&2\end{pmatrix},\quad
t\mapsto\begin{pmatrix}1&0\\0&\zeta_3\end{pmatrix},\quad
u\mapsto\frac{\zeta_3}{(\ii+1)\sqrt{3}}
\begin{pmatrix}-2\zeta_3&\zeta_3\\2&2\end{pmatrix}
\\\hline G_{12}&\vrule width0pt height 0.7cm depth 0.5cm 
s\mapsto D=\frac12
\begin{pmatrix}1&1+\frac1{\sqrt{-2}}\\2+\sqrt{-2}&-1\end{pmatrix},\quad
t\mapsto\overline D,\quad
u\mapsto\begin{pmatrix}1&0\\ 0&-1\end{pmatrix}
\\\hline G_{13}&\vrule width0pt height 0.7cm depth 0.5cm 
s\mapsto\begin{pmatrix}1&0\\0&-1\end{pmatrix},\quad
t\mapsto\frac1{\sqrt2}\begin{pmatrix}1&-1\\-1&-1\end{pmatrix},\quad
u\mapsto\frac1{\sqrt2}\begin{pmatrix}1&-\ii\\ \ii&-1\end{pmatrix}
\\\hline G_{14}&\vrule width0pt height 0.7cm depth 0.5cm 
s\mapsto\begin{pmatrix}1&0\\ 0&-1\end{pmatrix},\quad
t\mapsto\frac{\zeta_3^2}2
\begin{pmatrix}-1+\sqrt{-2}&1\\-1&-1-\sqrt{-2}\end{pmatrix}
\\\hline G_{15}&\vrule width0pt height 0.7cm depth 0.5cm 
s\mapsto C,\quad
t\mapsto\begin{pmatrix}1&0\\ 0&\zeta_3\end{pmatrix},\quad
u\mapsto\frac1{\sqrt3}\begin{pmatrix}1&\zeta_3^2\\2\zeta_3&-1\end{pmatrix}
\\\hline G_{16}&\vrule width0pt height 0.7cm depth 0.5cm 
s\mapsto\begin{pmatrix}1&0\\ 0&\zeta_5\end{pmatrix},\quad
t\mapsto\frac1{\sqrt5} 
\begin{pmatrix}1-\zeta_5^3&\zeta_5^4-1\\
\zeta_5-\zeta_5^2&\zeta_5-\zeta_5^3\end{pmatrix}
\\\hline G_{17}&\vrule width0pt height 0.7cm depth 0.5cm 
s\mapsto M=\frac{\ii}{\sqrt5}
\begin{pmatrix}\zeta_5^4-\zeta_5&\zeta_5^3-\zeta_5^2\\
\zeta_5^3-\zeta_5^2&\zeta_5-\zeta_5^4\end{pmatrix},\quad  
t\mapsto\begin{pmatrix}1&0\\ 0&\zeta_5\end{pmatrix}
\\\hline G_{18}&\vrule width0pt height 0.7cm depth 0.5cm 
s\mapsto N=\frac{\zeta_3^2}{\sqrt5}
\begin{pmatrix}\zeta_5^2-\zeta_5^4&1-\zeta_5^4\\
\zeta_5-1&\zeta_5^3-\zeta_5\end{pmatrix},\quad  
t\mapsto\begin{pmatrix}1&0\\ 0&\zeta_5\end{pmatrix}
\\\hline G_{19}&\vrule width0pt height 0.7cm depth 0.5cm 
s\mapsto N,\quad 
t\mapsto M,\quad
u\mapsto\begin{pmatrix}1&0\\0&\zeta_5\end{pmatrix}
\\\hline G_{20}&\vrule width0pt height 0.7cm depth 0.5cm 
s\mapsto\begin{pmatrix}1&0\\ 0&\zeta_3\end{pmatrix},\quad 
t\mapsto\frac{\zeta_3^2}{2\sqrt{-15}}
\begin{pmatrix}-5-\sqrt{-15}&2\\10&5-\sqrt{-15}\end{pmatrix}
\\\hline G_{21}&\vrule width0pt height 0.7cm depth 0.5cm 
s\mapsto\frac1{2\sqrt3}
\begin{pmatrix}1+\sqrt5&-1+\frac1{\sqrt5}\\-5+\sqrt5&-1-\sqrt5\end{pmatrix},\quad
t\mapsto\begin{pmatrix}1&0\\ 0&\zeta_3\end{pmatrix}
\\\hline
\end{longtable}

\subsection*{Example : the case of $G_9$}

Let us start with the model
$$ s\mapsto \begin{pmatrix}0&-\zeta_{24}^{11} \\ \zeta_{24} & 0\end{pmatrix},
\quad t\mapsto \begin{pmatrix}\ii& -\zeta_3 \\ 0 &  1\end{pmatrix}
$$
which is defined over $\BQ(\zeta_{24})$, an
extension of degree 2 of $K = \BQ(\zeta_8)$. By Galois descent or directly
one finds that conjugating by $\begin{pmatrix}1&0\\0&\zeta_3\end{pmatrix}$
yields the following model over $K$.
$$ s\mapsto\begin{pmatrix}0&-\zeta_8\\ \zeta_8^3 & 0\end{pmatrix},
\quad t\mapsto\begin{pmatrix} \ii & -1 \\ 0 &  1\end{pmatrix}
$$
We have $\Gamma \simeq \BZ/2 \times \BZ/2$, with three subgroups
of order 2 generated by $\gal{-1},\gal5$ and $\gal3$. To
these three subgroups correspond three intermediate extensions
between $\BQ$ and $K$. For each of these extensions, we get
a cocycle $\gal{k} \mapsto A_k \in GL_2(K)$, given by the following matrices
$$
A_{-1} = \begin{pmatrix} -\ii &  0 \\ 0 & 1 \end{pmatrix},\quad
A_5 = \begin{pmatrix} -1 & -1-\ii \\ 0 & 1 \end{pmatrix},\quad
A_3 = \begin{pmatrix} \ii & -1+\ii \\ 0 & 1 \end{pmatrix}
$$
It turns out that the last step of the algorithm is easy, because
these three cocycles match, so we get
in fact an element of $Z^1(\Gamma,\GL_2(K))$. In order to apply
Hilbert Theorem 90 we need to find $\lambda \in K$ such that
$\lambda + \gal{-1}(\lambda)A_1 + \gal5(\lambda)A_2 + \gal3(\lambda)A_3$
is invertible ;  it happens that $\lambda = 1 + \zeta_8 + \zeta_8^3$ is a solution
that gives the following invariant model
$$
s \mapsto \frac12 \begin{pmatrix}
-1-\zeta_8+\zeta_8^3 & 1/2-\zeta_8+\zeta_8^3 \\ 2 & 1+\zeta_8-\zeta_8^3  
\end{pmatrix},\quad
t \mapsto\begin{pmatrix} \ii & \frac{\ii -1 }2 \\ 0 & 1  \end{pmatrix}
$$
Finally, we know that the subgroup generated by $t$ is stable
by $\tilde{\iota}(\Gamma)$. It is thus possible to diagonalize it by a
rational matrix. We get the still invariant and simpler model
$$
s \mapsto\frac12 \begin{pmatrix} \sqrt2 & 2 \\ 1 & -\sqrt2  \end{pmatrix},\quad
t \mapsto \begin{pmatrix} 1 & 0 \\ 0 & \ii  \end{pmatrix}
$$

\subsection*{An example for $G_7$}
The model $\rho$ of the reflection representation of $G_7$ given by
$$s\mapsto\begin{pmatrix}1&0\\0&-1\end{pmatrix},\quad
t\mapsto \frac{\zeta_3^2}{(1+\ii)}\begin{pmatrix}-1&-1\\ \ii&-\ii\end{pmatrix}
,\quad
u\mapsto \frac{\zeta_3^2}{(1+\ii)}\begin{pmatrix}-1&\ii\\ -1&-\ii\end{pmatrix}
$$
is     globally     invariant,     giving    rise     to    an    homomorphism
$\eta:\Gal(\BQ(\zeta_{12})/\BQ)\to\Aut(G_7)$  given by  (specifying automorphisms,
as  in table  \ref{titilde},   by the  images of $s,t,u$)
$\cc\mapsto(s,u\inv,t\inv),\gal 7\mapsto(s,s^{us},u^{ts})$.
However, the automorphism $\eta(\cc)$ is not compatible with the action of the
complex  conjugation on other characters of  $G$, contrary to the homomorphism
$\itilde$  given in table \ref{titilde}, which  corresponds to the model given
in  table \ref{tablemodels}. The  existence of this  example is related to the
fact   that  $G_7$  is  the  only   exceptional  group  such  that  the  image
$\ibar(\Gamma)$ does not commute with $N$.

\subsection{Proof of \ref{C} for exceptional groups}

In this section we will deduce the existence of an
$\itilde$-equivariant model for
every  irreducible representation of $G$ (outside the exceptional set of
theorem \ref{A}) from the existence of such a model for
the reflection representation $\rho_0$, by applying lemma \ref{mfree}. We will
need the following proposition.

\begin{proposition}  \label{t1} If  $G$ is  of type  $G(d,1,r)$ then,  for all
distinct  $\rho_1, \rho_2 \in  \Irr(G)$, we have  $(\rho_1 \otimes \rho_0 \mid
\rho_2)  \leq 1$.  This also  holds for  exceptional groups, except for $G_{27},
G_{29},G_{34}$,  and  $G_{36}=E_7$.  For  these  4 exceptional groups, we have
$(\rho_1  \otimes \rho_0 \mid \rho_2) \leq  2$ for distinct $\rho_1,\rho_2 \in
\Irr(G)$.
\end{proposition}
\begin{proof}
The proof for the exceptional groups is a case-by-case computer check. If $G$ is
of  type $G(d,1,r)$  with $r  \geq 2$,  let $H$  denote its  natural parabolic
subgroup  of  type  $G(d,1,r-1)$,  and  $\rho_1  \in \Irr(G)$. Then $(\Ind_H^G
\trivialrep  \mid \rho_0) = (\trivialrep \mid  \Res_H \rho_0) = 1$ because $H$
is  a  maximal  parabolic  subgroup.  It  follows that $\rho_1 \otimes \rho_0$
embeds in $\rho_1 \otimes \Ind_H^G \trivialrep \simeq \Ind_H^G \Res_H \rho_1$.
In particular, for $\rho_2 \in \Irr(G)$ we have
$$
(\rho_1 \otimes \rho_0 \mid \rho_2) \leqslant (\Ind_H^G \Res_H \rho_1 \mid \rho_2) = (\Res_H \rho_1 \mid \Res_H \rho_2)
$$
so  it is  sufficient to  check that,  if $\rho_1  \neq \rho_2$, then $(\Res_H
\rho_1  \mid \Res_H \rho_2) \leq 1$. This is obvious considering the branching
rule (see \cite[p. 104]{Zelevinsky}).
\end{proof}

\subsection*{\label{rational  for primitive} Proof of \ref{C} for $G(d,1,r)$
and exceptional groups distinct from $G_{27},G_{29}$ or $G_{34}$}

We  already proved that there exists $\tilde{\iota} : \Gamma \to \Aut(G)$ such
that  
\begin{enumerate}
\item $\forall\rho\in\Irr(G)$, $\gamma\circ\rho\simeq\rho\circ\ibar(\gamma)$.
\item  $\rho_0$ admits a model  such that $\gamma \circ  \rho_0 = \rho_0 \circ
\itilde(\gamma)$.
\end{enumerate}
We  want to show that all $\rho  \in \Irr(G)$ admit an $\itilde$-equivariant model over
$K$,  that is  a model  over $K$  such that  $\gamma \circ  \rho =  \rho \circ
\tilde{\iota}(\gamma)$  for all $\gamma  \in \Gamma$. This  obviously holds if
$G$  is a Weyl group, so  we may assume that $G$  is not $G_{36} = E_7$. First
note that every  $\rho \in  \Irr(G)$ embeds  in some $\rho_0^{\otimes n}$ for
some  $n$, because $\rho_0$ is a faithful representation of $G$ (\cf.
\cite[problem 2.37]{FH}). It follows
that we can define the \emph{level} $N(\rho) \in \BN$ of $\rho$ by
$$
N(\rho) = \min \{ n \in \BN \ \mid \ \rho \into \rho_0^{\otimes n} \}.
$$
In  particular, $N(\rho) = m+1$ implies  that there exists $\rho' \in \Irr(G)$
with  $N(\rho') = m$ such that $\rho$ embeds in $\rho' \otimes \rho_0$. By (i)
we  know that  the representations  of level  at most  1, that  is the trivial
representation  and  $\rho_0$,  admit  an  equivariant  model.  We  proceed by
induction  on the level. Assume that all representations of level $m$ admit an
equivariant  model, and let $\rho \in \Irr(G)$  such that $N(\rho) = m+1$. Let
$\rho'  \in \Irr(G)$ of level $m$ such that $\rho \into \rho' \otimes \rho_0$.
By  proposition  \ref{t1}  we  have  $(\rho  \mid  \rho' \otimes \rho_0) = 1$.
Moreover  $\rho' \otimes \rho_0$ admits an  equivariant model over $K$ because
$\rho'$  and $\rho_0$ do. Then lemma \ref{mfree} implies that $\rho$ admits an
equivariant model over $K$ and we conclude by induction on the level.

\subsection*{Proof of \ref{C} for $G_{27},G_{29},G_{34}$ }

Let  $\CE \subset \Irr(G)$  denote the set  of exceptional (rational)
representations  described  in  the  introduction.  We note that $\trivialrep,
\rho_0 \not \in \CE$. Let us consider the following algorithm.

\medskip
\noindent {\bf Algorithm } 
\begin{enumerate}
\item  $L \leftarrow \{ \trivialrep, \rho_0 \}$.
\item For all $\rho' \in L$, $L' \leftarrow L \cup \{ \rho \in \Irr(G) \setminus \CE \ \mid \ (\rho' \otimes \rho_0 \mid \rho) = 1 \}$.
\item If $\# L' > \# L$, then $L \leftarrow L'$ and go to (ii).
\item Return $L$.
\end{enumerate}

The  same arguments as above, based on lemma \ref{mfree}, show that
all representations in the subset
$L$ of $\Irr(G) \setminus \CE$ returned by this algorithm admit an equivariant
model over $K$. Note that this algorithm only uses the character table of $G$.
It is then enough to check that it returns $\Irr(G) \setminus \CE$
for  $G_{27},G_{29}$  and  $G_{34}$  in  order  to  conclude  the proof of the
theorem.  Indeed, this  is the  case if  one applies  step (ii) five times for
$G_{27}$ and $G_{34}$, and four times for $G_{29}$.

\subsection{The exceptional case of $G_{22}$}
\label{sectionG22}

In this case we have $K = K' = \BQ(\ii,\sqrt5)$. We start from
the following model $\rho : G \hookrightarrow \GL(V)$ over
the degree 2 extension $K''=\BQ(\zeta_{20})$ of $K$.

$$
 s\mapsto\frac{\zeta_5^2-\zeta_5}{\zeta_{20}\sqrt 5}
\begin{pmatrix} 1 & \frac{1+\sqrt5}2 \\ \frac{1+\sqrt5}2 & -1 \end{pmatrix},
\quad t \mapsto \begin{pmatrix} 0 & \zeta_{20}^9 \\ -\zeta_{20} & 0
\end{pmatrix}, \quad 
u \mapsto \begin{pmatrix} 0 & -\zeta_{20} \\ \zeta_{20}^9 & 0
\end{pmatrix}
$$

This model is $\Gal(K''/\BQ)$-globally invariant.
Note that the projection $\Gal(K'' /\BQ) \to \Gal(K/\BQ)$
is not split. We check this by noting that $\Gal(K'' /\BQ)
\simeq \BZ/4\times \BZ/2$ admits only one subgroup
isomorphic to $\Gal(K/\BQ) \simeq \BZ/2\times \BZ/2$, namely
$\Gal(K''/\BQ(\sqrt5))$, and this subgroup
does not surject on $\Gal(K/\BQ)$, because $\sqrt5 \in K$.

As  predicted by the theorem  of Benard and Bessis,  we manage to find a model
$\rho_0$  over $K$ by Galois descent. It is  possible to get a model
over  $K$ such that $\gal{-1} \circ \rho = \rho \circ \tilde{\iota}(\gal{-1})$
\emph{or}  such that $\gal7 \circ \rho = \rho \circ \tilde{\iota}(\gal7)$, but
not  both. The  reason is  that it  is not  possible to  get a model such that
$\gal{13} \circ \rho = \rho \circ \tilde{\iota}(\gal{13})$.

\begin{lemma} $\rho$ does not admit any model
over $K$ which is globally invariant by $\Gal(K/\BQ(\ii))$. 
\end{lemma}

\begin{proof}
We  assume by  contradiction that  we are  given such  a model $\rho_1$. Since
$\rho_1$ is faithful there exists $a \in \Aut(G)$ of order 2 such that $\rho_1
\circ  a = \gal{13} \circ \rho_1$. We  will see in table \ref{bessisaut} that
$G$ has no non-trivial central automorphisms, thus
all  automorphisms of  $G$ preserve  the set  of reflections,  thus by theorem
\ref{refl} $\Out(G)$ acts faithfully and transitively on the set of reflection
representations.  It follows  that $a$  and $\itilde(\gal{13})$  have the same
image   in   $\Out(G)$,   namely   $\ibar(\gal{13})$.   In   particular  $a  =
\itilde(\gal{13})   \circ  \Ad(h)$  for   some  $h  \in   G$.  Exhausting  all
possibilities, we find using $\rho_0$ that the image of the cocycle defined by
$a$ in $H^1(\Gal(K/\BQ(\ii)), \PGL_2(K)) \simeq K^{\times} / N(K^{\times})$ is
always  the class of $1/3$ modulo $N(K^{\times})$. Using \MAGMA\ we check that
this class is non trivial, thus leading to a contradiction.
\end{proof}

We  now investigate  which representations  of $G$  admit a globally invariant
model  over $K$. This group has  18 irreducible representations, including two
1-dimensional  ones. The  four 2-dimensional  (reflection) representations are
deduced  from  each  other  through  Galois  action,  so  none of them admits a
globally  invariant model over $K$. In addition to these, there are 2 faithful
4-dimensional representations and 2 faithful 6-dimensional representations.

The center of $G_{22}$ is cyclic of order 4, generated by $Z = (stu)^5$.
We have $Z^2 = (ts)^5$. The other representations have for kernel
a non-zero subgroup of $Z(G)$. 

The 8 odd-dimensional representations admit $\itilde$-equivariant models
over $K$ by proposition \ref{odd dimension}, using the argument of
proposition \ref{oddinvariant}.

A computer check shows that every non-faithful
irreducible representation appear with multiplicity 1 in
some tensor product $\rho_1 \otimes \rho_2$, where $\rho_1$
and $\rho_2$ are odd-dimensional irreducible representations.
It follows that all these admit $\itilde$-equivariant models.

We now study in some detail the fate of the remaining
(faithful) irreducible representations. In particular,
we prove the following.

\begin{proposition}
Let $\rho$ be one of the faithful 4-dimensional
or 6-dimensional representations of $G_{22}$. Then $\rho$
admit globally invariant models over $K$, which induce
morphisms $j_{\rho} : \Gamma \to \Aut(G)$. These morphisms can be chosen
injective, however the induced morphisms $\overline{j_{\rho}}
: \Gamma \to \Out(G)$ are never injective.
\end{proposition}

In both dimensions, one of the representations considered
here satisfy $\rho(Z)=-\ii$. It is sufficient to consider these
to prove the proposition, since they are conjugated in pairs
by the Galois action. We let $N = \dim \rho \in \{ 4,6 \}$.

Since $\rho$ is faithful, the existence of a globally invariant
model for $\rho$ leads to a morphism $j : \Gamma \to \Aut(G)$
satisfying $\gamma \circ \rho \simeq \rho \circ j(\gamma)$ for
all $\gamma \in \Gamma$,
such that the associated cocycle $J = \{ \gamma \mapsto J_{\gamma} \}$
is cohomologically trivial.

We use the following procedure to check the proposition for an
arbitrary $j : \Gamma \to \Aut(G)$ satisfying
$\gamma \circ \rho \simeq \rho \circ j(\gamma)$ for all $\gamma \in \Gamma$.
Let $\gamma_0 = \gal{13}$ and note that $K^{\gamma_0} = \BQ(\ii)$. Recall
that $\gamma_0$ and $\cc$ generate $\Gamma = \Gal(\BQ(\ii,\sqrt{5})/\BQ)$.
Restricting $J$
to $\langle\gamma_0\rangle$ yields a class in $\BQ(\ii)^{\times}/NK^{\times}$
where $N(x) = x \gamma_0(x)$. If this class is non-zero then $J$
cannot be cohomologically trivial; otherwise there
exists $M \in \GL_N(K)$, which can be explicitly determined, such that
$J_{\gamma_0} = M\inv \gamma_0(M)$ in $\PGL_N(K)$. Then the cocycle $\gamma \mapsto
Q_{\gamma}  = M J_{\gamma} \gamma(M)\inv$ is cohomologous to $J$ and
$Q_{\gamma_0} = 1$. Since
$$
\gamma_0(Q_{\cc}) = Q_{\gamma_0} \gamma_0(Q_{\cc}) = Q_{\gamma_0 \cc}
= Q_{\cc \gamma_0} = Q_{\cc} \cc(Q_{\gamma_0}) = Q_{\cc}
$$
it follows that $\gamma_0(Q_{\cc}) = Q_{\cc}$, i.e. $Q_{\cc}
\in \PGL_N(K)^{\gamma_0} = \PGL_N(K^{\gamma_0})$, and that
$Q_{\gamma_0 \cc} = Q_{\cc}$. In other words, $Q$ belongs to the
image of the inflation map $Z^1(\Gal(K^{\gamma_0}/\BQ),\PGL_N(K^{\gamma_0}))
 \to Z^1(\Gamma,PGL_N(K))$.
Since the inflation map between Brauer groups is injective this
map is also injective and we are again reduced to a cyclic case,
hence to a norm equation solvable with \MAGMA.

Using \CHEVIE\ we get all such morphisms. The list obtained does not
depend on whether $\rho$ has dimension 4 or 6. If $N=4$,
the morphisms $j$ such that $J$ is cohomologous to 0 all
satisfy $\overline{j(\cc)} = \overline{j(\gal{7})} = \ibar(\gal{7})$
and $\overline{j(\gal{13})} = \ibar(1)$. Some of these
morphisms $j$ are injective.

In case $N=6$, the morphisms $j$ such that $J$ is cohomologically trivial
are all injective, and again send $\gal{13}$ to an inner automorphism.
Moreover $\overline{j(\cc)} \in \{ \ibar(\cc), \ibar(\gal{7}) \}$
and the two possibilities occur.

The choice of any morphism such that $J$ is cohomologically trivial
leads to a globally invariant model for the representations.
One of these morphisms, such that $J$ is cohomologically trivial for
both representations, has the following simple form :
$j(\gal{13}) = \Ad(t)$ and $j(\cc)$ sends $(s,t,u)$ to
$(u^{tu},t,s^{tst})$.

\begin{remark} \label{benardbessis} 
\end{remark}
Using the results proved here, one gets a simpler
proof of the theorem of Benard and Bessis. First
consider non-dihedral $G(de,e,r)$. Since all
irreducible representations of the groups $G(de,e,r)$ appear
as multiplicity 1 component of some representation of
$G(de,1,r)$, we are reduced by lemma \ref{mfree} to the case
of $G(d,1,r)$, then to its reflection representation
by proposition \ref{t1}, and we know that this representation is defined over
$K$ by the very definition
of these groups, or by \cite[7.1.1]{benson} which apply to all cases.
The case of the dihedral groups is classical (their representations
are mostly reflection representations).
Most exceptional groups are then dealt with by proposition \ref{t1}.
Finally, the algorithm used here for the few remaining ones completes
the proof, and is more efficient than the one used in \cite{bessis}.

\section{On the morphisms $\ibar_{\chi}$}

We first show how to derive corollary \ref{critere  rationalite}, stated
for arbitrary reflection groups, from theorem \ref{A}. It is an obvious consequence
of the following proposition.

\begin{proposition}\label{B}  Let $G$ be any  finite complex reflection group,
let  $K$ be its field of definition and let $\Gamma=\Gal(K/\BQ)$. Let $\CS$ be
the  set of irreducible  characters of $G$  which cannot be realized over
$\BQ$. Then for any $\chi\in\CS$ there exists an injection
$\Gamma\xrightarrow{\ibar_\chi}\Out(G)$  such  that  for any $\gamma\in\Gamma$
$\gamma(\chi) = \chi\circ  \ibar_\chi(\gamma)$.
\end{proposition}

\begin{proof}[Proof of proposition \ref{B}]
Proposition  \ref{B} is  an obvious consequence of \ref{A}  when
$G$  is irreducible.  For a  general reflection group $G$, we
have a decomposition into irreducible groups: $V=V_1 \oplus \ldots\oplus V_n$,
$G=G_1\times\ldots\times  G_n$  where  $G_i\subset\GL(V_i)$.  If  $K_i$ is the
field  of definition of $G_i$  then $K$ is the  subfield of $\BC$ generated by
$K_1,\ldots,K_n$,  whence if $\Gamma_i=\Gal(K_i/\BQ)$ we have natural quotient
morphisms     $\Gamma\to\Gamma_i$     such     that     the     product    map
$\Gamma\to\prod_i\Gamma_i$   is   injective.   We   deduce  an  injective  map
$\Gamma\to\Out(G)$ by composing $\Gamma\to\prod_i\Gamma_i$ with the individual
maps  $\Gamma_i\to\Out(G_i)$ deduced from  \ref{A} and then  with the natural
injection  $\Out(G_1)\times\ldots\Out(G_n)\to\Out(G)$.  Since  any irreducible
character $\chi$ of $G$ is a product of irreducible characters of the $G_i$,
proposition \ref{B} follows readily.
\end{proof}

We  now  study  further  the  structure  of  $\Out(G)$,  w.r.t.  the morphisms
$\ibar_{V}  :  \Gamma  \to  \Out(G)$.  We  now assume that $G$ is irreducible.
Recall  that, when $G  \neq \Sgot_6$, we  have $\Out(G) =  \Cbar .\Abar$, as a
direct  consequence of theorem \ref{auts}, and $\Abar = \Nbar \rtimes \Gamma$.
Moreover  it is clear that $\Nbar$ normalizes $\Cbar$, and centralizes $\Cbar$
when $\Nbar$ acts trivially on the linear characters.

The next lemma is the cornerstone of the interaction between $\ibar_V(\Gamma)$
and $\Cbar$.
\begin{lemma} \label{act lineaire}
Let $G$ be an irreducible complex reflection group. Then
$\chi(\ibar_V(\gamma)(g))  = \gamma(\chi(g))$ for all linear characters $\chi$
and all $\gamma \in \Gamma$.
\end{lemma}
\begin{proof}
The lemma results from the fact that $\ibar_{\chi}$ can be chosen 
to be the same for $\chi_V$ and for the linear characters. 
This is clear from theorem \ref{A} since either $\ibar_\chi$ can be chosen
independent of $\chi$, or we are in a $G(de,e,r)$, and the action of $\Nbar$
is trivial on both $V$ and the linear characters, since they are the 
restriction from $G(de,1,r)$ of a character which does not split.
\end{proof}

\begin{proposition} $\ibar_{V}(\Gamma) \subset \Abar$ commutes with $\Cbar$. In particular
$\Abar$ commutes with $\Cbar$ as soon as $\Nbar = \{ 1 \}$.
\end{proposition}
\begin{proof} 
By theorems \ref{xx} and \ref{C} there exists an extension
$K''$  of  $K$,  a  $\BQ$-form  of  $V\otimes_K  K''$  and  a map
$\iota:\Gal(K''/\BQ)\to\Aut(G)$   such  that  for   $g\in  G\subset\GL(V)$  and
$\gamma\in\Gal(K''/\BQ)$   we   have   $\gamma(g)=\iota(\gamma)(g)$, which
induce $\ibar_V :
\Gal(K/\BQ) \to \Out(G)$, except possibly for $G_{22}$. In the case of $G_{22}$, from the globally invariant
model over $K''$ we get a morphism $\iota : \Gal(K''/\BQ) \to
\Aut(G)$. The induced morphism $j$ to $\Out(G)$ factors through
$\Gal(K/\BQ)$, as $\gamma \circ \chi_V = \chi_V$ for $\gamma \in \Gal(K''/K)$,
and $\Nbar = 1$ in this case.
It then coincides with $\ibar_V : \Gal(K/\BQ) \to \Out(G)$, as $\chi_V\circ j(\gamma)
= \chi_V\circ \ibar_V(\gamma)$ for all $\gamma \in \Gamma$,
hence $j(\gamma) \ibar_V(\gamma)\inv \in \Nbar = \{ 1 \}$.

Let $a_{\chi} \in C$
a central automorphism defined by $a_{\chi}(g) = g \chi(g)^{-1}$.
It follows that $\iota(\gamma)\circ
\alpha_\chi\circ\iota(\gamma)\inv(g)=
\iota(\gamma)(\iota(\gamma)\inv(g)\chi(\iota(\gamma)\inv(g)))=
g \iota(\gamma)(\chi(\iota(\gamma)\inv(g)))$.
Now $\chi\circ\iota(\gamma)\inv = \chi \circ\ibar_V(\gamma)\inv = \gamma\inv \circ \chi$ by lemma \ref{act lineaire},
and $\iota(\gamma)(\chi(\iota(\gamma)\inv(g))) =
\gamma(\chi(\iota(\gamma)\inv(g)))$ as $\chi(\iota(\gamma)\inv(g)) \in ZG \subset G$.
It follows that $\iota(\gamma)(\chi(\iota(\gamma)\inv(g))) = \chi(g)$ hence
$\iota(\gamma)\circ
\alpha_\chi\circ\iota(\gamma)\inv(g) = \alpha_\chi(g)$ whence $\ibar_V(\Gamma)$ commutes to $C$.
\end{proof}


\section{Invariants}
We prove corollary \ref{invariants}, whose statement we recall.
We denote by $V^*$ the dual of a vector space $V$.
\begin{proposition} Let $G\subset\GL(V)$ be an irreducible complex reflection
group where $V$  is a $K''$-vector space, with $K''$
as in the introduction. 
There   is  a  $\BQ$-form  $V=V_0\otimes_\BQ  K''$  such  that  the  fundamental
invariants  of $G$  can be  taken rational,  \ie. in  the symmetric algebra
$S(V_0^*)$.
\end{proposition}
\begin{proof}
By  theorem \ref{C} and the matrix model over $K''$ exhibited for
$G_{22}$ in section \ref{sectionG22},
we may  assume that  $G\subset\GL(V_0\otimes_\BQ K'')$ is
globally  invariant  by  $\Gamma=\Gal(K''/\BQ)$.  Let  $f_1,\ldots,f_r$ (where
$r=\dim   V$)  be  fundamental   invariants,  \ie.  algebraically  independent
polynomials   such  that  $S(V^*)^G=K''[f_1,\ldots,f_n]$.   We  want  to  find
algebraically independent $g_1,\ldots,g_r\in S(V_0^*)$ such that we still have
$S(V^*)^G=K''[g_1,\ldots,g_n]$.  Our strategy will be  as follows: we will set
$g_i=\sum_{\gamma\in\Gamma}  \gamma(\lambda f_i)$ where $\lambda\in K''$. Then
the $g_i$ are still invariant since for $g\in G$ we have
$g(g_i)=\sum_{\gamma\in\Gamma}  \gamma(\lambda \gamma\inv(g)f_i) =\sum_{\gamma
\in\Gamma} \gamma(\lambda f_i)$, the last equality since $\gamma\inv(g)\in G$.

It is thus sufficient to show that we may choose $\lambda$ such that $g_i$ are
still  algebraically independent. By  the Jacobian criterion,  it is enough to
show  that we may choose $\lambda$ such that, if $x_1, \ldots, x_r$ is a basis
of  $V_0^*$, so that $S(V^*)\simeq  K''[x_1,\ldots,x_r]$, we have $\det\left (
\frac{\partial g_i}{\partial x_j} \right)_{i,j} \ne 0$.

We use the following version of the algebraic independence of automorphisms:
\begin{proposition}\cite[Chap. V, \S10, Th{\'e}or{\`e}me 4]{Bbk1}.
Let $Q$ be an infinite field, let $K$ be a finite Galois extension of $Q$ with
Galois  group $\Gamma$, let $\Omega$ be an  arbitrary extension of $K$ and let
$\{X_\gamma\}_{\gamma\in\Gamma}$  be indeterminates indexed by the elements of
$\Gamma$.  Let $F\in\Omega[X_\gamma]_{\gamma\in\Gamma}$  be a  polynomial such
that $F((\gamma(x))_{\gamma\in\Gamma})=0$ for any $x\in K$. Then $F=0$.
\end{proposition}
We   apply  the  proposition  with  $Q=\BQ$,  $\Omega=K''$  and  $F=\det\left  (
\sum_{\gamma\in\Gamma}   X_\gamma  \gamma(\frac{\partial   f_i}{\partial  x_j}
)\right)_{i,j}$.  The polynomial $F$ evaluated at $X_1=1$ and $X_\gamma=0$ for
$\gamma\ne  1$  is  equal  to  $\det\left  ( \frac{\partial f_i}{\partial x_j}
\right)_{i,j}$  which is non-zero,  so $F$ is  non-zero. By the theorem, there
exists  $\lambda\in  K''$  such  that  $F((\gamma(\lambda))_{\gamma})\ne 0$. But
$F((\gamma(\lambda))_{\gamma})=  \det\left ( \frac{\partial g_i}{\partial x_j}
\right)_{i,j}$ \end{proof}

We now prove \ref{quotient} whose statement we recall.
\begin{proposition} The morphism $V^\reg\to V^\reg/G$ is defined over $\BQ$.
\end{proposition}
\begin{proof}
We use the notations of the previous proof, in particular we 
assume that $G\subset\GL(V_0\otimes_\BQ K'')$ is
globally invariant by $\Gamma=\Gal(K''/\BQ)$.

Let $\CH$ be the set of reflecting hyperplanes for $G$ and for each $H\in\CH$,
let  $l_H$ be a linear form defining $H$. For $H\in\CH$ let $e_H$ be the order
of the subgroup of $G$ fixing $H$, and let $\Delta=\prod_{H\in\CH} l_H^{e_H}$.
It     is    well    known    that     $\Delta\in    S(V^*)^G$    (see    \eg.
\cite[6.44]{orlik-terao}).  Thus the variety $V^\reg/G$ is the open subvariety
of  $\Spec(K''[f_1,\ldots,f_n])$ whose function ring  is the localization by the
principal  ideal $\Delta$. It  is thus enough  to show that  we may choose the
$f_i$  such that a multiple of  $\Delta$ belongs to $\BQ[f_1,\ldots,f_n]$. If,
as  in the  previous proposition,  we choose  $f_i\in S(V_0^*)^G$,  it will be
enough  to  show  that  we  have  a  multiple  of $\Delta$ in $S(V_0^*)$ since
$S(V_0^*)\cap  S(V^*)^G=S(V_0^*)^G$.  Since  $\CH$  is  globally  invariant by
$\Gamma$,  for any  $\gamma\in\Gamma$ there  exists $\lambda_\gamma\in K''$ such
that   $\gamma(\Delta)=\lambda_\gamma   \Delta$,   and   it   is   clear  that
$\{\gamma\mapsto\lambda_\gamma\}\in   Z^1(\Gamma,  K^{\prime\prime\times})$.   By  Hilbert's
theorem   90   we   have   that   $H^1(\Gamma,K^{\prime\prime\times})$   is   trivial,  thus
$\{\gamma\mapsto\lambda_\gamma\}$   is   a   coboundary,   \ie.  there  exists
$\lambda\in  K^{\prime\prime\times}$ such  that $\lambda_\gamma=\lambda\inv\gamma(\lambda)$.
We   get  then  that   $\gamma(\lambda\inv\Delta)=\lambda\inv\Delta$  for  any
$\gamma$, thus $\lambda\inv\Delta\in S(V_0^*)$.
\end{proof}
\section{Quasi-indecomposable groups, and non-irreducible groups}

In this chapter, we prove  theorems \ref{quasi-indec} and \ref{KRS for
quasi-indec}, and then go on to describe the non-abelian factors and the
central automorphisms of irreducible complex reflection groups.

We recall that a group is indecomposable if it has no non-trivial decomposition
as a direct product. We recall the statement of the Krull-Remak-Schmidt
theorem for finite groups.

\begin{theorem}[Krull-Remak-Schmidt]\label{KRS}
Given  two  decompositions  $G=G_1\times\ldots\times G_k=H_1\times\ldots\times
H_{k'}$  of the finite  group $G$ as  a product of  indecomposable factors, we
have  $k=k'$ and  there is  a permutation  $\sigma$ of  $\{1,\ldots,k\}$ and a
central automorphism $\alpha$ such that $\alpha(H_i)= G_{\sigma(i)}$.
\end{theorem}
\begin{proof}
See \eg.  \cite[3.3.8]{robinson}.
\end{proof}

\subsection{Quasi-indecomposable groups}
We  say that  a group  $G$ is  quasi-indecomposable if  in any  direct product
decomposition $G=G_1\times G_2$, we have either $G_1\subseteq ZG$ or
$G_2\subseteq ZG$.

\begin{proposition}\label{quasidec}
Let $G$ be a finite quasi-indecomposable group. Then $G$ admits a
decomposition of the form $G=Z\times \hat{G}$ where $Z\subseteq ZG$ and $\hat{G}$ is
either trivial or is a non-abelian indecomposable group. In such a
decomposition $\hat{G}$ (called the {\em non-abelian factor} of $G$ and denoted
$\can(G)$) and $Z$
(called the {\em central factor} of $G$) are unique up to isomorphism.
\end{proposition}
\begin{proof}
First,  we observe that $G$ admits such a decomposition. Indeed, being finite,
$G$  admits a  decomposition $G=G_1\times\ldots  \times G_n$  as a  product of
indecomposable groups. By assumption, at most one of the $G_i$ is non-abelian;
we  set $\hat{G}$ equal  to that factor  and $Z$ equal  to the product of the other
factors.

Assume  now that  $G$ admits  two such  decompositions $G=Z\times \hat{G}=Z_1\times
\hat{G}_1$. If we refine the decompositions to a product of indecomposable groups,
by the Krull-Remak-Schmidt theorem there is an automorphism which sends one
decomposition to the other, and this automorphism must send the non-abelian
factor in one decomposition to the other, whence the isomorphisms between
$\hat{G}$ and $\hat{G}_1$ and $Z$ and $Z_1$.
\end{proof}

We now relate decompositions of quasi-indecomposable groups to central
endomorphisms. An endomorphism $\alpha:
G\to G$ is {\em central} if for any $g \in G$, we have $\alpha(g)\in g ZG$.
\begin{lemma}
Let  $G$ be a finite quasi-indecomposable group, and let $\alpha$ be a central
endomorphism  of $G$ with a maximal kernel (that is, $\ker\alpha$ is a maximal
element  for inclusion among the $\ker\beta$ where  $\beta$ runs over the
central   endomorphisms).   Then   $G=\ker\alpha\times\Image(\alpha)$  is  a
decomposition of $G$ with $\ker\alpha \subseteq ZG$  and $\Image(\alpha)$
indecomposable non abelian (or trivial).
\end{lemma}
\begin{proof}
Since  $\ker\alpha\subseteq\ker\alpha^2$  and  $\ker\alpha$  is maximal, we have
$\ker\alpha^2=\ker\alpha$,   whence   $\Image(\alpha^2)=\Image(\alpha)$  since
$\#G=\#\ker\alpha  \#\Image(\alpha)$. In  particular $\alpha$  is injective on
$\Image(\alpha)$,     so    $\ker\alpha\cap     \Image(\alpha)=\{1\}$.    Thus
$G=\ker\alpha\times\Image(\alpha)$;  and  $\Image(\alpha)$  is indecomposable,
otherwise  if $\Image(\alpha)=Z\times \hat{G}$ where $Z$ is central, and if $\beta$
is the projection $\Image(\alpha)\to \hat{G}$, then $\beta\circ\alpha$
would be
central and would have a
greater kernel than $\alpha$, contradicting the maximality of $\ker\alpha$.
\end{proof}

Given a  central endomorphism  of $G$,  we define an homomorphism
$\chi_\alpha:G\to   ZG$   by   the   formula  $\chi_\alpha(g)=g\alpha(g\inv)$.
Conversely  any  homomorphism  $\chi:G\to  ZG$  gives rise to a central
endomorphism   by   the   formula   $\alpha_\chi(g)=g\chi(g)\inv$.  We  have
$\ker\alpha_\chi=\{z\in ZG\mid \chi(z)=z\}$.

\begin{lemma} Let $G$ be a quasi-indecomposable finite group with cyclic
center. Then in a decomposition $G=Z\times \hat{G}$ with $Z$ central and $\hat{G}$
indecomposable nonabelian, the factor $Z$ is unique.
\end{lemma}
\begin{proof}
In view of the previous lemma, it is sufficient to prove that given two
central endomorphisms $\alpha$ and $\beta$, there is a central automorphism
$\gamma$ such that $\ker\alpha\subseteq\ker\gamma$ and
$\ker\beta\subseteq\ker\gamma$. Let us take $\gamma=\alpha\circ\beta$. Then
clearly $\ker\beta\subseteq\ker\gamma$, and we will have 
$\ker\alpha\subseteq\ker\gamma$ if we can show that restricted to $ZG$ the
endomorphisms $\alpha$ and $\beta$ commute.
Now $\alpha(\beta(z))=\beta(z)\chi_\alpha(\beta(z))\inv=
z\chi_\beta(z\inv)\chi_\alpha(z\inv)\chi_\alpha(\chi_\beta(z))$, and this 
formula is symmetric in $\alpha$ and $\beta$ since the two endomorphisms 
$\chi_\alpha$ and $\chi_\beta$ of the
cyclic group $ZG$ commute with each other.
\end{proof}
\begin{remark} Non-unicity of the central factor.
\end{remark}
The group $\BZ/2\times W(B_2)$
does not have a unique central factor. 
We will see that $W(B_2)$ is indecomposable, so the above is a 
decomposition into a central and  a non-abelian indecomposable factor.
If $\{1,\varepsilon\}$ are the two
elements of the factor $\BZ/2$ and $\{1,w_0\}$ the two elements of
$Z(W(B_2))$ then $\{1,\varepsilon w_0\}$ is another direct factor
of $W(B_2)$.
\begin{remark} Non-unicity of the non-abelian factor.
\end{remark}
Even if the central factor is unique, $\can(G)$ may not be.
We will see that $W(G_2)$ is quasi-indecomposable. Having a cyclic center,
it has a unique central factor which is the center $\{1,w_0\}$.
A possible non-abelian factor is $W(A_2)\subset W(G_2)$. However, this is not
the only one. If we consider the sign character $x \mapsto \sgn(x)$ of $W(A_2)$ to take
its values in $\{1,w_0\}$ then the subgroup formed of $\{x\sgn(x)\mid
x\in W(A_2)\}$ is another possible non-abelian factor.

We will see however that both factors are unique for most irreducible
complex reflection groups.
We prove now proposition \ref{quasi-indec}.
\begin{proposition}
Irreducible finite complex reflection groups are quasi-indecomposable.
\end{proposition}
\begin{proof}
We  prove  the  result  case  by  case,  starting  with the imprimitive groups
$G(de,e,r)$  of  rank  $r\ge  3$.  Recall  that  such  a  group is of the form
$G=D\rtimes\Sgot_r$,  where $D$ is  the subgroup of  diagonal matrices in $G$.
Let  $\pi$  be  the  quotient  map  $G\to\Sgot_r$ and consider a decomposition
$G=G_1\times  G_2$. Then we  get $\Sgot_r=\pi(G_1)\pi(G_2)$, a  product of two
normal  subgroups.  Since  any  proper  normal  subgroup  of  $\Sgot_r$  is in
$\Agot_r$,  we must  have \eg.  $\pi(G_1)=\Sgot_r$. Since  $\pi(G_2)$ commutes
with  $\pi(G_1)$ and the  center of $\Sgot_r$  is trivial (since  $r\ge 3$) we
have   $\pi(G_2)=1$,   \ie.   $G_2\subseteq   D$.   Further,   the   action   of
$\pi(G_1)=\Sgot_r$ on $G_2$ is trivial, so $G_2\subseteq ZG$.

To  treat the  other cases,  we will  follow the  arguments of \cite[6.2 lemme
2]{marin}.   Assume that $G\subset\GL(V)$,   an   irreducible   group,  admits  a
decomposition  $G_1\times G_2$ where neither of the $G_i$ is abelian. Since an
irreducible  representation of a product is  the tensor product of irreducible
representations  of the factors,  we must have  $V=V_1\otimes V_2$ where $V_i$
affords a faithful irreducible representation of $G_i$. In particular, we must
have $\dim V_i>1$ since $G_i$ is non-abelian.

These remarks immediately solve the cases of $G(de,e,2)$ or more generally the
irreducible  groups  of  prime  rank.  It  remains to consider the exceptional
groups  of  rank  $4,6$  or  $8$.  For  that,  we  use  that  if  $G$  is  not
quasi-indecomposable  then $V\otimes V\simeq  V_1\otimes V_1\otimes V_2\otimes
V_2=(S^2  V_1\oplus\Lambda^2 V_1)\otimes (S^2  V_2\oplus\Lambda^2 V_2)$ has at
least  4  irreducible  components  (since  $\dim  V_i>  1$  implies  that  the
decomposition $V_i\otimes V_i=S^2 V_i\oplus\Lambda^2 V_i$ is non-trivial).

Now  it  is  a  theorem  of  Steinberg  \cite[Exercice  3]{Bbk2}  that  for an
irreducible  complex reflection group, $\Lambda^2  V$ is irreducible. Thus, we
see  that if  $S^2V$ has  less than  3 irreducible  components, then $V\otimes
V=S^2  V\oplus  \Lambda^2  V$  will  have  less  than  4, and thus $G$ will be
quasi-indecomposable.  When $G$ is real (a  Coxeter group) the
representation $S^2V$ contains the trivial representation.  It turns out that for {\em any}
exceptional  irreducible complex reflection group (not  only those of rank 4,6
or  8),  $S^2V$  is  irreducible  when  $G$  is  not  real  and  $S^2V-\trivialrep$ is
irreducible  when $G$ is  real (this can  easily checked by  computer from the
character table of $G$).
\end{proof}

We have the following consequence of the Krull-Remak-Schmidt theorem
for direct products of quasi-indecomposable groups:
\begin{theorem}\label{KRS qindec} 
Let   $G=G_1\times\ldots\times   G_n$   be   a   direct   product   of  finite
quasi-indecomposable    groups,   such   that   for    any   $i,j$   we   have
$\can(G_i)\simeq\can(G_j)\Rightarrow  G_i\simeq G_j$. Then any automorphism of
$G$   is  of  the  form  $\alpha\circ\beta$,   where  $\alpha$  is  a  central
automorphism   and  $\beta$  satisfies   $\beta(G_i)=G_{\sigma(i)}$  for  some
permutation $\sigma$ of $\{1,\ldots,n\}$.
\end{theorem}
\begin{proof}
Let $\gamma \in \Aut(G)$ and let $G=G_1\times\ldots\times G_k$ be a decomposition of
$G$ in indecomposable factors. Let us apply the Krull-Remak-Schmidt theorem to
the  decomposition given by  $H_i=\gamma(G_i)$; we get  a central automorphism
$\alpha$  of  $G$  and  $\sigma  \in  \mathfrak{S}_k$ such that $\alpha(G_i) =
H_{\sigma\inv(i)}$.  Letting $\beta = \alpha\inv \circ \gamma$ it follows that
$\gamma$  can  be  written  $\alpha  \circ  \beta$  with  $\alpha$  a  central
automorphism  of $G$  and $\beta(G_i)  = G_{\sigma(i)}$  for some  $\sigma \in
\mathfrak{S}_k$.

If we lump together the product of all abelian $G_i$ as a single subgroup $Z$,
we  get that  any automorphism  of a  group of  the form $Z \times G_1 \times
\dots  \times  G_n$  where  $Z$  is  abelian  and  the  $G_i$  are non-abelian
indecomposable  groups is of  the form $\alpha  \circ \beta$ where $\alpha$ is
central and $\beta(G_i) = G_{\sigma(i)}$ for some $\sigma \in \mathfrak{S}_n$,
as   $\beta(G_i)$,  being  non-abelian,  cannot  be   mapped  to  one  of  the
indecomposable  factors of  $Z$. Let  $\beta_Z \in  \Aut(G)$ being  defined by
$\beta_Z(g)  = \beta(g)$ for $g \in Z$ and  $\beta_Z(g) = g$ for $g \in G_i$ ;
similarly,  $\beta_G(g) = \beta(g)$ for  $g \in G_i$, $\beta_G(g)  = g$ for $g
\in  Z$. We have $\beta  = \beta_Z \circ \beta_G$  and $\gamma = (\alpha \circ
\beta_Z)  \circ \beta_G$. Since $\beta_Z$ is a central automorphism, replacing
$\alpha$ by $\alpha \circ \beta_Z$ and $\beta$ by $\beta_G$ we thus can assume
that $\beta$ acts trivially on $Z$.

Rewriting  the  group  $G$  of  the  statement  as  $Z_1\times\ldots\times Z_n\times
\hat{G}_1\times  \ldots\times \hat{G}_n$ where $Z_i$ are the  central factors of the $G_i$ and
$\hat{G}_i$  are their non-abelian factors, we get that any automorphism of $G$ is of
the  form $\alpha\circ\beta$ where $\alpha$ is  central and $\beta$ is trivial
on  $Z_1\times\ldots\times  Z_n$  and  effects  a  permutation $\sigma$ of the
$\hat{G}_i$.       The       assumption       of       the      statement      that
$\can(G_i)\simeq\can(G_j)\Rightarrow   G_i\simeq   G_j$   now   implies   that
$G_{\sigma(i)}\simeq G_i$, hence by \ref{quasidec} there exists isomorphisms 
$\delta_i : Z_{\sigma(i)} \to  Z_i$ between the corresponding central
factors.
We recall that $\beta$ acts
trivially  on  each  $Z_i$,  and  introduce  $\delta  \in  \Aut(G)$  such that
$\delta(g)  = g$  for $g  \in \hat{G}_i$  and $\delta(g)  = \delta_i(g)$ for $g \in
Z_i$. Letting $\beta' = \delta\inv \circ \beta$ we have $\beta'(g) = \beta(g)
\in  \hat{G}_{\sigma(i)}$ for  $g \in  \hat{G}_i$ and  $\beta'(g) = \delta_i\inv(g) \in
Z_{\sigma(i)}$  for  $g  \in  G_i$.  Since  $\delta$ is clearly central we get
$\gamma  =  \alpha'  \circ  \beta'$,  with  $\alpha'  = \alpha \circ \delta$ a
central  automorphism and  $\beta'(G_i) =  G_{\sigma(i)}$, which concludes the
proof.
\end{proof}
We now prove theorem \ref{KRS for quasi-indec}. 

\begin{theorem}
Let  $G\subset\GL(V)$ be a complex reflection group and let $G=G_1\times\ldots
\times G_n$  be its  decomposition in  irreducible factors. 
Then the following are equivalent:
\begin{enumerate}
\item No $G_i$ is isomorphic to $\Sgot_6$, and for any $i,j$, we
have $\can(G_i)\simeq\can(G_j)\Rightarrow G_i\simeq G_j$.
\item For any $\rho\in\Irr(G)$, any automorphism of $G$ is of the form 
$\alpha\circ\beta\circ\gamma$, where $\alpha$ is central, $\beta$ is induced
by $N_{\GL(V)}(G)$, and $\gamma$ stabilizes the $G_i$ and
preserve their reflections. 
\end{enumerate}
Moreover, for every complex reflection group $G$, any automorphism
which preserve the reflections
can be written as $\beta \circ \gamma$, with $\beta,\gamma$ as above.
\end{theorem}
\begin{proof}
We recall that $N_{\GL(V)}(G)$, mapping reflections to reflections,
necessarily permutes the factors $G_i$. We first observe that $(ii)\Rightarrow(i)$ is clear, since if one of the $G_i$
is isomorphic to $\Sgot_6$ the exceptional automorphisms of $\Sgot_6$ are not of
the  form given in (ii), and, if two non-isomorphic components $G_i$ and $G_j$
have  isomorphic non-abelian factors: $G_i=Z_i\times \hat{G}$ and $G_j=Z_j\times \hat{G}$,
the  automorphism of  $G$ which  exchanges the  isomorphic factors  $\hat{G}$ while
leaving the other factors invariant is not either of the form given in (ii).


By \ref{KRS qindec} any automorphism is of
the  form  $\alpha_1\circ\beta_1$  where  $\alpha_1$  is central and $\beta_1$
permutes the $G_i$. Since a permutation of isomorphic $G_i$ can be effected by
an   element  of   $N$,  we   can  write  $\beta_1=\beta_2\circ\beta_3$  where
$\beta_2\in N$ and $\beta_3(G_i)=G_i$. This proves the last assertion of the
theorem. Moreover, since no $G_i$ is isomorphic to $\Sgot_6$,
by \ref{auts}, $\beta_3$ itself is
of  the  form  $\beta_3=\alpha_2\circ\gamma$,  where $\alpha_2$ is
central and  $\gamma$ preserves the reflections. 
Since the subgroup $C$
of    central    automorphisms    is    normal,       $\alpha_1\circ\beta_2
\circ\alpha_2\circ\gamma$ is of the required form. This concludes
the proof of the theorem.

\end{proof}
\section{Central endomorphisms, central automorphisms and non-abelian factors
of irreducible complex reflection groups}
In  this section, we describe the central factor and the group $C$ of
central automorphisms for irreducible complex reflection groups.

We  relate  the  problem,  for  complex  reflection  groups,  to  linear
characters.  
We recall that the group of linear characters, or equivalently the abelianized
of   $G$,  is  isomorphic  to  $\prod_H  \BZ/\#C_G(H)$  where  $H$  runs  over
representative  of hyperplane orbits. If $S$  is the generating subset we have
taken for $G$, the hyperplane orbits define a partition $\sigma$ of $S$, and a
linear    character   $\chi$   is    specified   by   specifying   arbitrarily
$\chi(s)\in\mu_{e_s}$ for one $s$ in each part of $\sigma$, where $e_s$ is the
order of $s$.

Such  a  character  $\chi$  will  give  rise  to  a  central  endomorphism  if
$\chi(s)\in\mu_{\#ZG}$  for  any  $s\in  S$.  Let  $z$ be a generator of $ZG$;
$\alpha_\chi$ will be an automorphism if and only if $z\chi(z)\inv$ is still a
generator of $ZG$. Notice that $\ker \alpha_{\chi} \subseteq ZG$.

We  now proceed, case-by-case, to describe $\sigma$ and  give an expression
for $z$ in terms of the generators $S$, from which we will deduce $C$ and the
central factor of $G$ as well as a possible non-abelian factor.

For  the  groups  $G(de,e,r)$,  we  have $S=\{t',s'_1,s_1,\ldots,s_{r-1}\}$ (where
$t'=t^e$  and $s'_1=s_1^t$ as in section  \ref{Gde,e,r}) where $t'$ is omitted
if  $d=1$ and $s'_1$ is omitted if  $e=1$. Expressions for a generator of $ZG$
can  be found  in \cite{BMR};  it is  equal to  $z=(t's_1\ldots s_{r-1})^r$ if
$e=1$,        and        $z=t^{\prime        r/\gcd(e,r)}       (s'_1s_1\ldots
s_{r-1})^{\frac{e(r-1)}{\gcd(e,r)}}$  in the general  case (this expression is
still   valid  for   $d=1$  by   setting  $t'=1$   in  that   case).  We  have
$ZG\simeq\mu_{d\cdot\gcd(e,r)}$.

When  writing the partition  $\sigma$, we decorate  the parts where $e_s\ne 2$
with  the  value  of  $e_s$,  separated  by  a  semicolon.  Thus  $\sigma$  is
$\{t';d\},\{s'_1,s_1,\ldots,s_{r-1}\}$  except when $r=2$ and $e$ is
even, in which case it is 
$\{t';d\},\{s'_1\},\{s_1\}$.
We  find  that  a  character  $\chi$  gives  rise to a central endomorphism of
$G(de,e,r)$ if and only if $\chi(s'_1)=\chi(s_i)=1$ if $d\gcd(e,r)$ is odd.
We assume this condition in the sequel. We have
$\chi(z)=\chi(t')^{\frac    r{\gcd(e,r)}}$. 

We first consider the case where $r>2$  or $r=2$ and  $e$ is odd.
If    we    define    $j$    by
$\chi(t')=\zeta_d^j$ and $0 \leq j < d$,  and  if  $\alpha_\chi$  is  the  corresponding  central
endomorphism, we then have that
\begin{multline*}
\#\ker(\alpha_\chi)=\#\{i\in\{0,\ldots,d\gcd(e,r)-1\}\mid
\chi(t')^{\frac{ri}{\gcd(e,r)}}=\zeta_{d\gcd(e,r)}^i\}=\\ \hfill
\#\{i\in\{0,\ldots,d\gcd(e,r)-1\}\mid i(rj-1)\equiv 0\pmod{d\gcd(e,r)}\}\\
\end{multline*}

We define $d_{r'}$ to be the largest factor of $d$ which is prime to
$r$; notice that $rj-1$ and $d \gcd(e,r)$ being coprime is
equivalent to $rj-1$ and $d_{r'}$ being coprime.

We thus find that
\begin{itemize}
\item  $\alpha_\chi$ is  an automorphism  if and  only if  $rj-1$ is  prime to
$d_{r'}$.
\item The maximal kernel for an $\alpha_\chi$ is of cardinality $d_{r'}$,
obtained when $rj\equiv 1\pmod {d_{r'}}$. There always exists such a
$\chi$ as $d_{r'}$ and $r$ are coprime.
\end{itemize}

If  $rj\equiv  1\pmod  {d_{r'}}$  we  have  $\det(\alpha_\chi(t'))  =\det(t')/
\chi(t')^r=    \zeta_d^{1-rj}\in\mu_{d/d_{r'}}$,   thus
$\Image(\alpha_\chi) \subseteq G(de,d_{r'} e,r)$ unless it is possible
to take  $\chi(s_i)=-1$ and $-1\notin G(de,d_{r'}e,r)$, which happens when
$d$ is even and $r$ odd. We finally get that 
when $r>2$ or $r=2$ and $e$ is odd the central factor of
$G(de,e,r)$ is $\BZ/d_{r'}$ and the non-abelian factor is unique equal to
$G(de,d_{r'}e,r)$, unless $d$ is even and $r$ odd when there is another
possibility for the non-abelian factor.

In  the remaining  case of  $G(2de',2e',2)$, the partition  $\sigma$ is
$\{t';d\},\{s'_1\},\{s_1\}$, and $z=t'(s'_1s_1)^{e'}$ is of order $2d$.

Define  $\epsilon=0$ if $\chi(s'_1  s_1)=1$ or $e'$  is even, and $\epsilon=d$
otherwise,  and define $j$ by $\chi(t')=\zeta_d^j$. An analysis similar to the
above   case   shows   that   $\#\ker  \alpha_\chi=\#\{i\in\{0,\ldots,2d-1\}\mid
i(2j+\epsilon-1)\equiv0\pmod{2d}\}$.

We  find that $\chi$ defines a central  automorphism if and only if $\epsilon$
is  even  and  $2j+\epsilon-1$  is  prime  to  $d_\odd$, where $d_\odd$ is the
largest odd divisor of $d$.

The  maximal kernel  of $\alpha_\chi$  is of  order $d_{2'}$ obtained \eg. for
$\epsilon=0$ and $2j\equiv 1\pmod {d_{2'}}$, except when $d$ and $e'$ are both
odd,  in which case it is of  order $2d$, obtained when $\chi(s'_1s_1)=-1$ and
$j\equiv  (1+d)/2\pmod d$. It can be checked  that in this last case there
are two possibilities for the image of $\alpha_\chi$; if
$\chi(s'_1)=-1$ it lies in $G(de',de',2)$, so we get
\begin{proposition}
Decompositions of $G(de,e,r)$ into a central factor and a
non-abelian indecomposable factor  are as follows:
\begin{itemize}
\item
If  $r=2$,  $d$  is  odd  and  $e=2e'$  with  $e'$  odd,  a  decomposition  is
$\BZ/2d\times   G(de',de',2)$;  there   is  one   other  possibility  for  the
non-abelian factor.
\item    In   the   other   cases   a   decomposition   is   $\BZ/d_{r'}\times
G(de,d_{r'}e,r)$; this decomposition is unique unless $d$ is even and $r$ odd,
in which case there is one other possibility for the non-abelian factor.
\end{itemize}
\end{proposition}

Table  \ref{bessisaut}  gives  for  exceptional  groups the partition $\sigma$
(where  each part is decorated by  the order of the corresponding reflections,
when this order is not 2), the order of $ZG$, an expression in terms of $S$ of
a  generator of $ZG$ (such expressions can be found in \cite[table 4]{BMR} and
in   \cite{bessis-michel}),  under  the   column  ``automorphism''  conditions
describing the linear characters giving rise to an element of $C$ (this column
is empty when there are no conditions), and under the column ``decomposition''
the decomposition of $G$. In this table, $\mathfrak{A}_5.2$
denotes the (unique) nonsplit central extension of $\mathfrak{A}_5$
by $\BZ/2$, and $\SO_5(\BF_3)'$ is the commutator subgroup of
$\SO_5(\BF_3)$.

This decomposition turns out to be unique in each
case, excepted for $G_5$ and $G_7$ where the non-abelian factor given is one 
of 3 possibilities.
\begin{longtable}[c]{|>{$}l<{$}|>{$}c<{$}>{$}c<{$}>{$}c<{$}>{$}c<{$}>{$}c<{$}|}
\caption{\label{bessisaut}Central endomorphisms}\\
\hline 
G&\sigma&\#ZG&z&\text{automorphism}&\text{decomposition}\\
\hline 
\endfirsthead
\hline 
G&\sigma&\#ZG&z&\text{automorphism}&\text{decomposition}\\
\hline 
\endhead
\hline 
\endfoot
G_4&\{s,t;3\}&2&(st)^3&\chi(s)=1&G_4\\
G_5&\{s;3\},\{t;3\}&6&(st)^2&\chi(st)\ne\zeta_3&\BZ/3\times G_4\\
G_6&\{s\},\{t;3\}&4&(st)^3&\chi(t)=1&G_6\\
G_7&\{s\},\{t;3\},\{u;3\}&12&stu&\chi(tu)\ne\zeta_3&\BZ/3\times G_6\\
G_8&\{s,t;4\}&4&(st)^3&&G_8\\
G_9&\{s\},\{t;4\}&8&(st)^3&&G_9\\
G_{10}&\{s;3\},\{t;4\}&12&(st)^2&\chi(s)\ne\zeta_3&\BZ/3\times G_8\\
G_{11}&\{s\},\{t;3\},\{u;4\}&24&stu&\chi(t)\ne\zeta_3&\BZ/3\times G_9\\
G_{12}&\{s,t,u\}&2&(stu)^4&&G_{12}\\
G_{13}&\{s\},\{t,u\}&4&(stu)^3&&G_{13}\\
G_{14}&\{s\},\{t;3\}&6&(st)^4&\chi(t)\ne\zeta_3&\BZ/3\times G_{12}\\
G_{15}&\{s\},\{t\},\{u;3\}&12&u(st)^2&\chi(t)\ne\zeta_3&\BZ/3\times G_{13}\\
G_{16}&\{s,t;5\}&10&(st)^3&\chi(s)\ne\zeta_5^2&\BZ/5\times(\Agot_5.2)\\
G_{17}&\{s\},\{t;5\}&20&(st)^3&\chi(t)\ne\zeta_5^2&\BZ/5\times G_{22}\\
G_{18}&\{s;3\},\{t;5\}&30&(st)^2&\chi(s)\ne\zeta_3,\chi(t)\ne\zeta_5^2&\BZ/15\times
(\Agot_5.2)\\
G_{19}&\{s\},\{t;3\},\{u;5\}&60&stu&\chi(t)\ne\zeta_3,\chi(u)\ne\zeta_5^2&
\BZ/15\times G_{22}\\
G_{20}&\{s,t;3\}&6&(st)^5&\chi(s)\ne\zeta_3&\BZ/3\times(\Agot_5.2)\\
G_{21}&\{s\},\{t;3\}&12&(st)^5&\chi(t)\ne\zeta_3&\BZ/3\times G_{22}\\
G_{22}&\{s,t,u\}&4&(stu)^5&&G_{22}\\
G_{23}&\{s,t,u\}&2&(stu)^5&\chi(s)=1&\BZ/2\times\Agot_5\\
G_{24}&\{s,t,u\}&2&(stu)^7&\chi(s)=1&\BZ/2\times\GL_3(\BF_2)\\
G_{25}&\{s,t,u;3\}&3&(stu)^4&&G_{25}\\
G_{26}&\{s\},\{t,u;3\}&6&(stu)^3&\chi(s)=1&\BZ/2\times G_{25}\\
G_{27}&\{s,t,u\}&6&(uts)^5&\chi(s)=1&\BZ/2\times\Agot_6\\
G_{28}&\{s,t\},\{u,v\}&2&(stuv)^6&&G_{28}\\
G_{29}&\{s,t,u,v\}&4&(stvu)^5&&G_{29}\\
G_{30}&\{s,t,u,v\}&2&(stuv)^{15}&&G_{30}\\
G_{31}&\{s,t,u,v,w\}&4&(stuvw)^6&&G_{31}\\
G_{32}&\{s,t,u,v;3\}&6&(stuv)^5&\chi(s)\ne\zeta_3^2&\BZ/3\times\Sp_4(\BF_3)\\
G_{33}&\{s,t,u,v,w\}&2&(stvwu)^9&\chi(s)=1&\BZ/2\times \SO_5(\BF_3)'\\
G_{34}&\{s,t,u,v,w,x\}&6&(stvwux)^7&&G_{34}\\
G_{35}&\{s_1,\ldots,s_6\}&1&1&\chi(s)=1&G_{35}\\
G_{36}&\{s_1,\ldots,s_7\}&2&(s_1\ldots
s_7)^9&\chi(s)=1&\BZ/2\times\SO_7(\BF_2)\\
G_{37}&\{s_1,\ldots,s_8\}&2&(s_1\ldots s_8)^{15}&&G_{37}\\
\end{longtable}

We remark that for any irreducible complex reflection group, the non-abelian
factor of $G$ can always be taken to be a reflection subgroup, unless $G$ is an exceptional
group in the list
$G_{16},G_{18},G_{20},G_{23},G_{24},G_{27},G_{32},G_{33},G_{36}$.

The   information  we  give  is  sufficient,  in  rank  2,  to  determine  the
intersection  of $C$  with the  automorphisms which  preserve the reflections,
thus  determining the  structure of  $\Aut(G)$ in  that case.  For exceptional
groups,  the isomorphism type of $C$ is  given in table 1 of \cite{bessis} (we
note  on this table  that the group  structure of $C$  given by composition of
automorphisms    is   unrelated to the   group   structure   on
$\Hom(G,\BC^\times)$;  in particular $C$ may be non-commutative).

Assume  that  $G$  is  of  rank  $2$  and  let  $H$  be  the  group  of Galois
automorphisms,   which   is   well-determined   by   theorem   \ref{xx}  since
$\gcd(e,r)\le  2$ in our case. It follows  from the above table that $C=H$ for
$G_{12},G_{13},G_{14},G_{15}$,   that   $H$   is   of   index  2  in  $C$  for
$G_8,G_9,G_{10},G_{11}$,  that $C$ is  of index 2  in $H$ for $G_4,G_6,G_{20},
G_{20},G_{22}$,  and that $C\cap  H$ is of  index $2$ in  both $C$ and $H$ for
$G_{16},G_{17},G_{18},G_{19}$. Finally, in the two cases $G_5$ and $G_7$ where
$N$  is not trivial, we have $N\subset C$. In these
cases   $C$  is  non-commutative,  isomorphic   to  $\Sgot_3$  for  $G_5$  and
$\BZ/2\times\Sgot_3$  for $G_7$. In both cases $C\cap H$ is the center of $C$,
which is trivial for $G_5$ and equal to $\BZ/2$ for $G_7$.





\end{document}